\documentclass[12pt]{article}
\usepackage[ansinew]{inputenc}
\usepackage[pdftex]{graphicx}
\usepackage{amssymb,amsmath,amscd}
\usepackage{dsfont}
\usepackage{verbatim}
\usepackage{amsmath} 
\usepackage{graphicx}
\usepackage{pstricks}
\usepackage{enumitem}
\usepackage{ntheorem}
\usepackage{amssymb}
\usepackage{float}
\usepackage{breqn}
\usepackage[left=2cm,top=2cm,right=2cm,bottom=2cm,nohead]{geometry}
\newtheorem{de}{Definition}[section]

\newtheorem{theo}{Theorem}
  
\newtheorem{lem}{Lemma}[section]
\newtheorem{cor}[lem]{Corollary} 
\newtheorem{conj}{Conjecture}[section]

\setenumerate[1]{font=\bfseries,label=\Roman*.}
\setenumerate[2]{font=\itshape,label=\arabic*)}

\author{Yashar Memarian
}       

\title{A Geometric Approach to Radial Correlation Type Problems \\
       }
\date{\today}

\begin{document}
\maketitle
\tableofcontents

\begin{abstract}
A radial probability measure is a probability measure with a density (with respect to the Lebesgue measure) which depends only on the distances to the origin. Consider the Euclidean space enhanced with a radial probability measure. A correlation problem concerns showing whether the radial measure of the intersection of two symmetric convex bodies is greater than the product of the radial measures of the two convex bodies. A radial measure satisfying this property is said to satisfy the \emph{correlation property}. A major question in this field is about the correlation property of the (standard) Gaussian measure. The main result in this paper is a theorem suggesting a sufficient condition for a radial measure to satisfy the correlation property. A consequence of the main theorem will be a proof of the correlation property of the Gaussian measure.

\end{abstract}

\section{Introduction}

Consider the Euclidean space $\mathbb{R}^n$ enhanced with a probability measure $\mu_n$. Assume this probability measure has a density function which is a radial function (depending only on the distances to the origin). We write down this probability measure as $\mu_n=g(r)dx$ where $dx$ is the Lebsegue measure and $g(r)$ is a (continuous) function depending only on the distances to the origin. We call such a measure a \emph{radial probability measure}. By definition, a symmetric convex body $K$ is a convex set in $\mathbb{R}^n$ which contains the origin and which is symmetric with respect to it. The general question of interest is:

\emph{Question}:

For which radial probability measure $\mu_n$ the following property holds: for any two symmetric convex bodies $K_1,K_2\subset \mathbb{R}^n$, we have:
\begin{equation} \label{eqn:ay}
\mu_n(K_1\cap K_2)\geq \mu_n(K_1)\mu_n(K_2).
\end{equation}

A probability measure which satisfies inequality (\ref{eqn:ay}) (for every pair of symmetric convex bodies) is said to satisfy the \emph{correlation property}.

The main radial probability measure for which the above question was intensively studied is the \emph{Gaussian measure} where the standard Gaussian measure (denoted by $\gamma_n$) of any measurable subset $A\subseteq \mathbb{R}^n$ is defined by
\begin{equation} \label{eqn:gaus}
\gamma_n(A)=\frac{1}{(2\pi)^{n/2}}\int_{A}e^{-\vert x\vert^2/2} dx.
\end{equation}
The Gaussian Correlation conjecture was a captivating problem in the field of convex geometry and probability theory. A \emph{less} general form of the Gaussian Correlation Conjecture first appeared in $1955$ in \cite{dunnet} and \emph{more} general form appeared a few years after in $1972$ by S. Das Gupta, M.L. Eaton, I. Olkin, M. Perlman, L.J. Savage and M. Sobel, in \cite{dasgu}. The one-dimensional case of the conjecture was proven by both Khatri and Sidak (independently) in $1966-67$ in \cite{khatri}, \cite{sidak} and \cite{sidakk}. The two-dimensional case was proven by Pitt in $1977$ in \cite{pitt}. Harg\'e proved the conjecture for ellipsoids in $1999$ (see \cite{harge}). Schechtman, Schlumprecht and Zinn proved partial results in \cite{sh}. The subject was quite active and several specialists have worked on variations of the conjecture. At last, Thomas Royen in \cite{roy} successfully presented a complete proof of the conjecture. Impressively, Royen's proof is a very short proof.

The main theorem of this paper is the following:

\begin{theo} \label{strong}
Let $n\geq 4$. Let $\mu_n$ be a radial probability measure where $\mu_n=g(r)dx$ and $dx$ is the Lebesgue measure. If for every two symmetric strips $S_1$ and $S_2$ in $\mathbb{R}^2$ where the axis of either $S_1$ or $S_2$ coincides with the $x$-axis we have:
\begin{equation} \label{eqn:do}
\mu_2(\mathbb{R}^2)\mu_2(S_1\cap S_2)\geq \mu_2(S_1)\mu_2(S_2),
\end{equation}
where $\mu_2=\vert y\vert^{n-2}g(r)dx\,dy$, then for every two symmetric convex bodies $K_1$ and $K_2$ in $\mathbb{R}^n$, we have :
\begin{eqnarray*} 
\mu_n(K_1\cap K_2)\geq \mu_n(K_1)\mu_n(K_2),
\end{eqnarray*}
Here, the measure $\mu_2$ is defined by
\begin{eqnarray*}
\mu_2=\vert y\vert^{n-2}g(r)dx\,dy,
\end{eqnarray*}
and a symmetric strip $S\subset \mathbb{R}^2$ is an open, convex, symmetric set such that a $u\in \mathbb{S}^1$ and $h>0$ exist such that
\begin{eqnarray*}
S=\{x\in\mathbb{R}^2: \vert x.u\vert<h\}.
\end{eqnarray*} 
The axis of a strip $S$ is the axis directed by the unit vector $u\in\mathbb{S}^1$.
\end{theo}

Theorem \ref{strong} gives a sufficient condition for a radial probability measure to satisfy the correlation property. This condition is a verification of a $2$-dimensional analytic problem.

The idea for the proof of the above theorem is by using localisation on the Riemannian sphere, and more precisely, bringing the $n$-dimensional inequality (\ref{eqn:ay}) down to a $2$-dimensional correlation problem for symmetric $2$-dimensional convex sets (with respect to the anisotropic measure presented in the statement of Theorem \ref{strong}). This will indeed simplify the high dimensional complexity of the problem and leaves one with a $2$-dimensional geometric problem. By generalising the geometric ideas presented recently in the paper \cite{figa}, one can go further and simplify the geometry of symmetric convex sets to only consider linear geometric objects, which are the symmetric strips. This again will simplify the $2$-dimensional geometry and leaves us with a linear analytic problem.

This paper is organised as follows: In section $3$, an overview of the idea of the proof of theorem \ref{strong} will be presented. The aim would be to give the main idea of the proof before delving into the details and the techniques of the next sections. In Section $4$, the theory of convexly-derived measures and the localisation on the canonical Riemannian sphere are recalled. The tehchniques and results presented in this section are the one enabling one to simplify the study of theorem \ref{strong} to the study of a family of $2$-dimensional correlation problems. Section $5$ concerns the study of a family of $2$-dimensional correlation type problems. The main idea of this section is the extension of the symmetisation procedure used in \cite{figa} to a family of $2$-dimensional anisotropic correlation problems.  In Section $6$, the proof of the main theorem \ref{strong} will be completed. In Section $7$, as a corollary of theorem \ref{strong}, the proof of the Gaussian Correlation Problem will be presented. Finally, in Section $8$, one explains how the methods used in this paper can be applied for other problems in convex geometry most notably the Mahler Conjecture.

\section{Acknowledgement}
It should be mentioned that I first learned of this conjecture via Olivier Gu\'edon, whom I thank. I am also grateful to M.Ledoux who pointed out that I ought to look at \cite{figa}. I am very grateful to K.Grove who motivated me to edit and resubmit this paper.

\section{Idea of the Proof of Theorem \ref{strong}}

It may be useful to give a summary of the way Theorem \ref{strong} will be proved. The details of everything written in this section can be found in the following chapters. 

The space upon which we are working is the Euclidean space $\mathbb{R}^n$. The metric on this space is supposed to be the canonical  Euclidean metric. The measure, however, is a radial probability measure, which can be expressed as $\mu_n=g(r) dx$ where $dx$ is the Lebsegue measure in $\mathbb{R}^n$ and $g(r)$ is a continuous function which only depends on the distances to the origin.
The geometric objects on which we will be computing their measures are convex bodies which contain the origin and are symmetric with repect to it. 

The question of interest is the following : given any pair of symmetric convex bodies $K_1$ and $K_2$ in $\mathbb{R}^n$, when does the following hold?
\begin{equation} \label{eqn:ti}
\mu_n(\mathbb{R}^n)\mu_n(K_2\cap K_2)\geq \mu_n(K_1)\mu_n(K_2).
\end{equation}
Of course, $\mu_n(\mathbb{R}^n)=1$ but for reasons which will become clear shortly, the expression appears in the left hand side of inequality (\ref{eqn:ti}). The inequality (\ref{eqn:ti}) is an inequality involving four continuous functions defined on the canonical sphere $\mathbb{S}^{n-1}$. Indeed it is sufficient to write down the integrals in polar coordinates associated to $\mathbb{R}^n$:
\begin{eqnarray*} 
\mu_n=r^{n-1}g(r)dr\,du_{\mathbb{S}^{n-1}},
\end{eqnarray*}
where $du_{\mathbb{S}^{n-1}}$ is the canonical Riemannian measure associated to the canonical sphere with sectional curvature everywhere equal to $1$. Inequality (\ref{eqn:ti}) is simply the following:
\begin{equation} \label{eqn:ti2}
\int_{\mathbb{S}^{n-1}}f_1(u)du \int_{\mathbb{S}^{n-1}}f_2(u)du\geq \int_{\mathbb{S}^{n-1}}f_3(u)du \int_{\mathbb{S}^{n-1}}f_4(u) du,
\end{equation}
where for $i=1,2,3,4$, $f_i$ is a continuous function defined on $\mathbb{S}^{n-1}$ defined as follows: set $M_1=\mathbb{R}^n$, $M_2=K_1\cap K_2$, $M_3=K_1$ and $M_4=K_2$. The function $f_i$ takes a point $u\in\mathbb{S}^{n-1}$ and calculates the following :
\begin{eqnarray*}
f_i(u)=\int_{0}^{x_i(u)}r^{n-1}g(r)dr,
\end{eqnarray*}
where $x_i(u)$ is the the length of the segment starting from $0$, going in the direction $u$, and touching the boundary of $M_i$. Of course for $M_1$, for every $u\in\mathbb{S}^{n-1}$, we have $x_1(u)=+\infty$. 

%So far, I've convinced the reader that inequality (\ref{eqn:ti}) is an inequality involving the integral of four continuous functions on $\mathbb{S}^{n-1}$. So, what next?

The next idea, which is the fundamental idea of this paper, is what is called \emph{the localisation technique}. The localisation technique basically states that if one wants to demonstrate an inequality such as inequality (\ref{eqn:ti2}), it is enough to prove it on \emph{one}-dimensional subsets of $\mathbb{S}^{n-1}$ (but with respect to a wider class of measures).

It turns out that the $1$-dimensional subsets are the geodesic segments of $\mathbb{S}^{n-1}$ which we denote one by $\sigma$. And it turns out that  the probability measures defined on $\sigma$ are written as $C\cos(t+t_0)^{n-2}dt$. Here we parametrise $\sigma$ by its arc length. The measure $dt$ is the canonical length measure of $\sigma$. The segment $\sigma$ is an interval with length at most equal to $\pi$ and on this interval, the measure has a density which is given by $\cos(t+t_0)^{n-2}$.

To sum up, applying the localisation technique, if one wants to prove inequality (\ref{eqn:ti2}), it is enough to prove the following inequality for every pair $(\sigma,\nu)$:
\begin{equation} \label{eqn:ti3}
\int_{\sigma}f_1(t)d\nu(t)\int_{\sigma}f_2(t)d\nu(t)\geq \int_{\sigma}f_3(t)d\nu(t)\int_{\sigma}f_4(t)d\nu(t).
\end{equation}
The functions $f_i$ were defined earlier and inequality (\ref{eqn:ti3}) concerns the restriction of these functions on $\sigma$. Taking the restriction of the functions $f_i$ on $\sigma$ means that one should intersects the sets $\mathbb{R}^n$, $K_1$, $K_2$ and $K_1\cap K_2$ with a $2$-dimensional plane which contains the origin of $\mathbb{R}^n$ and the geodesic segment $\sigma$. On this $2$-dimensional plane, take the cone which is defined over $\sigma$. This cone is denoted by $C$. Now one writes down the inequality (\ref{eqn:ti3}) on $C$. Parametrising $C$ in polar coordinates, one obtains the measure defined on $C$ which are given by $\mu_{2,t_0}=r^{n-1}f(r)\cos(t+t_0)^{n-2}dr\,dt$. Therefore, inequality (\ref{eqn:ti3}) is equivalent to the following:
\begin{equation} \label{eqn:ti4}
\mu_{2,t_0}(C\cap \mathbb{R}^2)\mu_{2,t_0}(K_1\cap K_2\cap C)\geq \mu_{2,t_0}(K_1\cap C)\mu_{2,t_0}(K_2\cap C).
\end{equation}

%We are then able to simplify the verification of Conjecture \ref{1} to the verifications of a family of $2$-dimensional problems given by inequality (\ref{eqn:ti4}). This means that if we can prove that for every cone $C$ and every measure $\mu_2$ of the form given above, if inequality (\ref{eqn:ti4}) is true then Conjecture \ref{1} is true.

Here, the $n$-dimensional problem is translated to a family of $2$-dimensional problems, but these $2$-dimensional problems can also have annoying geometry for which studying an inequality such as (\ref{eqn:ti4}) can become tricky.

% How do we study them then?

Next idea is to use an appropriate symmetrisation method for studying problems such as the one given by inequality (\ref{eqn:ti4}). 

It will be proved that for every $C$, $K_1$, $K_2$ and $\mu_{2,t_0}$ there are always two symmetric strips $S_1$ and $S_2$ such that:
\begin{equation} \label{eqn:axis}
\frac{\mu_{2,t_0}(C\cap \mathbb{R}^2)\mu_{2,t_0}(K_1\cap K_2\cap C)}{\mu_{2,t_0}(K_1\cap C)\mu_{2,t_0}(K_2\cap C)}\geq \frac{\mu_{2,t_0}(\mathbb{R}^2)\mu_{2,t_0}(S_1\cap S_2)}{\mu_{2,t_0}(S_1)\mu_{2,t_0}(S_2)}.
\end{equation}

Now if for every pair of symmetric strips $S_1$ and $S_2$  and every measure $\mu_{2,t_0}$:
\begin{equation} \label{eqn:ti5}
\mu_{2,t_0}(\mathbb{R}^2)\mu_{2,t_0}(S_1\cap S_2)\geq \mu_{2,t_0}(S_1)\mu_{2,t_0}(S_2),
\end{equation}
then inequality (\ref{eqn:ay}) holds true.

However it turns out that verifying equation (\ref{eqn:ti5}) for every pair of strips is too much information to consider. 

The last idea is to use the symmetry of the convex sets to show that it is enough to consider pair of strips where we have a control on the axis of at least one of them. This enormously simplifies the study of this problem. This part is where the parameter $t_0$ defined in the family of measures $\mu_{2,t_0}$ becomes important.

One shall show that for every $C$ there is always a $t_0$ such that inequality (\ref{eqn:axis}) holds, moreover, the axis of at least one strip coincides with the $x$-axis.

Therefore if for every pair of strips where the axis of at least one coincides with the $x$-axis, the following inequality:
\begin{equation} \label{eqn:ti6}
\mu_{2,t_0}(C\cap \mathbb{R}^2)\mu_{2,t_0}(K_1\cap K_2\cap C)\geq \mu_{2,t_0}(K_1\cap C)\mu_{2,t_0}(K_2\cap C),
\end{equation}
is true then inequality (\ref{eqn:ay}) is true and the proof of the theorem \ref{strong} follows.

%The good news is that it turns out that $C$ induced with these good measures $\mu_{2,t_0}$ are enough information for us to consider. This is how it works :

%Let us assume that Conjecture \ref{1} is not true. This means inequality (\ref{eqn:ti2}) is not true. In the next section, I will show that if (\ref{eqn:ti2}) is not true, then there exists a \emph{partition} of $\mathbb{S}^{n-1}$ into geodesic segments, such that on each of these geodesic segments a measure such as $\nu$ is defined, and for every pair of $(\sigma,\nu)$ in this partition the inequalities (\ref{eqn:ti3}), and hence (\ref{eqn:ti4}), are not true. The idea here is the following : if we cleverly obtain the partition of $\mathbb{S}^{n-1}$ into geodesic segments upon which inequalities (\ref{eqn:ti3}) and (\ref{eqn:ti4}) are not true, we can always find at least one geodesic segment $\sigma$ and a measure $\nu$ defined on this $\sigma$, such that translating the inequality on the cone $C$ defined over this $\sigma$, provides us with a measure $\mu_{2,t_0}$ (for which we know that inequality (\ref{eqn:ti6}) is \emph{always} true). This gives us a contradiction and finalises the demonstration of Conjecture \ref{1}.    

\subsection{Condensed Scheme of Proof}
The proof of Theorem \ref{strong} follows directly by applying the following three results which will be proved in the following sections: Lemma \ref{genlova}, Corollary \ref{corrr} and Lemma \ref{close}.
\begin{itemize}
\item Lemma \ref{genlova} is the localisation Lemma which enables one to bring the dimension $n$ down to $2$.
\item Corollary \ref{corrr} is the $2$-dimensional anisotropic correlation inequality which simplifies the study to the case of a pair of symmetric strips.
\item Lemma \ref{close} is the lemma which enables us to have a control on the axis of the symmetric strips proposed in Theorem \ref{strong}.
%\item Theorems \ref{twodim} and \ref{con} will tell us that for studying inequalities such as (\ref{eqn:ti4}) it is necessary to study similar inequalities involving only symmetric strips in $\mathbb{R}^2$.
%\item Lemma \ref{rot} and Corollary \ref{corrr} tell us that there is always measure $\mu_{2,t_0}$ for which the inequality (\ref{eqn:ti5}) is equivalent to an inequality for symmetric strips with a control on the axis of one of the strips. Here I call them, a pair of good strips.
%\item Lemma \ref{close} tells us that we can always find a partition of $\mathbb{S}^{n-1}$ into geodesic segments for which we assume inequality (\ref{eqn:ti4}) is not true, and find a pair of geodesic segments and a measure $\nu$ in this partition for which the inequality (\ref{eqn:ti5}) is as close as we want to the inequality involved with a pair of good strips.
%\item Lemma \ref{end} tell us that for every pair of good strips, the inequality (\ref{eqn:ti5}) is true.
%\item This gives us a contradiction and the proof of Theorem \ref{1} follows.

\end{itemize} 

One is ready to delve into the details and the different techniques used for the demonstration of theorem \ref{strong}. 

\section{Localisation on the Sphere} \label{sin}

In the past few years, localisation methods have been used to prove several very interesting geometric inequalities. In \cite{lova} and \cite{kann}, the authors proved integral formulae using localisation, and applied their methods to conclude a few isoperimetric-type inequalities concerning convex sets in the Euclidean space. In \cite{guedon} the authors study a functional analysis version of the localisation, used again on the Euclidean space. Localisation on more general spaces was studied in \cite{gromil}, \cite{grwst}, \cite{memwst}, and \cite{memusphere}.

Many materials in this section are derived from \cite{memwst}:

\begin{de}[Convexly-derived measures]
A convexly-derived measure on $\mathbb{S}^n$ (resp. $\mathbb{R}^n$) is a limit of a vaguely converging sequence of probability measures of the form $\mu_i=\frac{vol|S_i}{vol(S_i)}$, where $S_i$ are open convex sets. The space $\mathcal{MC}^n$ is defined to be the set of probability measures on $\mathbb{S}^n$ which are of the form $\mu_{S}=\frac{vol_{\vert S}}{vol(S)}$ where $S\subset \mathbb{S}^n$ is open and convex. The space of \emph{convexly-derived probability measures} on $\mathbb{S}^n$ is the closure of $\mathcal{MC}^n$ with respect to the vague (or weak by compacity of $\mathbb{S}^n$)-topology. The space $\mathcal{MC}^k$ will be the space of convexly-derived probability measures whose support has dimension $k$ and $\mathcal{MC}^{\leq k}=\cup_{l=0}^{k}\mathcal{MC}^{l}$. 
\end{de}

This class of measures was defined first in \cite{gromil} and used later on in \cite{ale}, \cite{memwst}, \cite{memphd}. In Euclidean spaces, a convexly-derived measure is simply a probability measure supported on a convex set which has a $x^k$-concave density function with respect to the Lebesgue measure. To understand convexly-derived measures on the sphere we will need some definitions:
\begin{de}[$\sin$-concave functions]
A real function $f$ (defined on an interval of length less than $2\pi$) is called $\sin$-concave, if, when transported by a unit speed paramatrisation of the unit circle, it can be extended to a $1$-homogeneous and concave function on a convex cone of $\mathbb{R}^2$.
\end{de}

\begin{de}[$\sin^k$-affine functions and measures]
A function $f$ is affinely $\sin^k$-concave if $f(x)=A\sin^k(x+x_0)$ for a $A>0$ and $0\leq x_0\leq \pi/2$. A $\sin^k$-affine measure by definition is a measure with a $\sin^k$-affine density function.
\end{de}

\begin{de}[$\sin^k$-concave functions]
A non-negative real function $f$  is called $\sin^k$-concave if the function $f^{\frac{1}{k}}$ is $\sin$-concave.
\end{de}

One can easily confirm the following:

\begin{lem} \label{sinconc}
\begin{itemize}
\item A real non-negative function$f$ defined on an interval of length less than $\pi$ is $\sin^k$-concave if for every $0<\alpha<1$ and for all $x_1,x_2 \in I$ we have
\begin{eqnarray*}
f^{1/k}(\alpha x_1+(1-\alpha)x_2)\geq (\frac{\sin(\alpha\vert x_2-x_1\vert)}{\sin(\vert x_2-x_1\vert)})f(x_1)^{1/k}+(\frac{\sin((1-\alpha)\vert x_2-x_1\vert)}{\sin(\vert x_2-x_1\vert)})f(x_2)^{1/k}.
\end{eqnarray*}
Particularly if $\alpha=\frac{1}{2}$ we have
\begin{eqnarray*}
f^{1/k}(\frac{x_1+x_2}{2})\geq \frac{f^{1/k}(x_1)+f^{1/k}(x_2)}{2\cos(\frac{\vert x_2-x_1\vert}{2})}.
\end{eqnarray*}
\item $f$ admits only one maximum point and does not have any local minima.
\item If $f$ is $\sin$-concave and defined on an interval containing $0$, then $g(t)=f(\vert t\vert)$ is also $\sin$-concave.
\item Let $0<\varepsilon<\pi/2$. Let $\tau >\varepsilon$. $f$ is defined on $[0,\tau]$ and attains its maximum at $0$. Let $h(t)=c\cos(t)^k$ where $c$ is choosen such that $f(\varepsilon)=h(\varepsilon)$. Then
\begin{eqnarray*}
\begin{cases}
f(x) \geq h(x)     & \text{for } x\in [0,\varepsilon], \\
f(x) \leq h(x)     & \text{for } x\in [\varepsilon,\tau].
\end{cases}
\end{eqnarray*}
In particular, $\tau\leq \pi/2$.
\item Let $\tau>0$ and $f$ be a nonzero non-negative $\sin^k$-concave function on $[0,\tau]$ which attains its maximum at $0$. Then $\tau\leq \pi/2$ and for all $\alpha\geq 0$ and $\varepsilon\leq \pi/2$ we have
\begin{eqnarray*}
\frac{\displaystyle\int_{0}^{min\{\varepsilon,\tau\}}f(t)dt}{\displaystyle\int_{0}^{\tau}f(t)dt}\geq \frac{\displaystyle\int_{0}^{\varepsilon}\cos(t)^{k}dt}{\displaystyle\int_{0}^{\pi/2}\cos(t)^{k}dt}.
\end{eqnarray*}
\end{itemize}
\end{lem}

This class of measures are also used in Optimal Transport Theory (see the excellent book \cite{villani} on this matter as well as a proof for Lemma \ref{sinconc}).

\begin{lem} \label{memlem}
Let $S$ be a geodesically convex set of dimension $k$ of the sphere $\mathbb{S}^n$ with $k \leq n$. Let $\mu$ be a convexly-derived measure defined on $S$ (with respect to the normalised Riemannian measure on the sphere). Then $\mu$ is a probability measure having a continuous density $f$ with respect to the canonical Riemannian measure on $\mathbb{S}^k$ restricted to $S$. Furthermore, the function $f$ is $\sin^{n-k}$-concave on every geodesic arc contained in $S$.
\end{lem}
The above Lemma, proved in \cite{memwst}, completely characterises the class of convexly-derived measures on the sphere. Note the similarity between the Euclidean case and the spherical one.

\begin{de}[Spherical Needles] \label{sn}
A spherical needle in $\mathbb{S}^n$ is a couple $(I,\nu)$ where $I$ is a geodesic segment in $\mathbb{S}^n$ and $\nu$ is a probability measure \emph{supported} on $I$ which has a $\sin^{n-1}$-affine density function.
\end{de}

\emph{Remark}:

According to definition \ref{sn}, one can properly writes down the measure $\nu$. To do so, choose a parametrisation of the geodesic segment $I$ by its arc length. Therefore there is a (canonical) map $s:[0,l(I)]\to I$. For every $t\in[0,l(s)]$, we have that $\Vert \frac{ds}{dt}\Vert=1$. The measure $dt$ is the canonical Riemannian length-measure associated to the geodesic segment $I$. Then, $I$ is parametrised by $t\in [0,l(I)]$, the measure $\nu$ can be written as $\nu=C\cos(t-t_0)^{n-1}dt$, for $t_0\in[0,\pi]$ and $C$ is the normalisation constant such that :
\begin{eqnarray*}
\int_{0}^{l(I)}C\cos(t-t_0)^{n-1}dt=1.
\end{eqnarray*}

It is necessary to say a few words on \emph{convex partitions}. The reason being the fact that later on, one needs the canonical sphere (seen as a metric measure space) to be partitioned into spherical needles. These objects need to be properly defined.
  
\begin{de} \label{convexparti}
Let $\Pi$ be a finite convex partition of $\mathbb{S}^n$. We review this partition as an atomic probability measure $m(\Pi)$ on the space $\mathcal{MC}$ as follows: for each piece $S$ of $\Pi$, let $\mu_{S}=\frac{vol_{\vert S}}{vol(S)}$ be the normalised volume of $S$. Then set 
\begin{eqnarray*}
m(\Pi)=\sum_{S}\frac{vol(S)}{vol(\mathbb{S}^n)}\delta_{\mu_{S}}.
\end{eqnarray*}
Define the \emph{space} of (infinite) \emph{convex partitions} $\mathcal{CP}$ as the vague closure of the image of the map $m$ in the space $\mathcal{P}(\mathcal{MC})$ of probability measures on the space of convexly-derived measures. The subset $\mathcal{CP}^{\leq k}$ of convex partitions of dimension $\leq k$ consists of elements of $\mathcal{CP}$ which are supported on the subset $\mathcal{MC}^{\leq k}$ of convexly derived measures with support of dimension (at most) $k$. It is worth remembering that the space $\mathcal{CP}$ is compact and $\mathcal{CP}^{\leq k}$ is closed within. 
\end{de}

%We now include few useful deﬁnitions and remarks which we shall need  later on. The point being, the element of the partition of the sphere will always consist of convex sets, and there exist an algorithmic procedure which enables us to construct such sets. The following deﬁnitions are related to this fact

\emph{Remark}:

There exists an algorithmic procedure which enables one to construct the elements of the partition. The following definitions are related to this fact.

\begin{de}[Pancakes]
Let $S$ be an open convex subset of $\mathbb{S}^n$. Let $\varepsilon>0$. We call $S$ an $(k,\varepsilon)$-pancake if there exists a convex set $S_{\pi}$ of dimension $k$ such that every point of $S$ is at distance at most $\varepsilon$ from $S_{\pi}$.
\end{de}

\emph{Remark}: 

Every geodesic segment $I$ is a Hausdorff limit of a sequence $\{S_i\}_{i=1}^{\infty}$ , where $S_i$ is a $(1,\varepsilon_i)$-pancake and where $\varepsilon_i\to 0$ when $i\to\infty$. Furthermore, every spherical needle $(I,\nu)$ is a limit of a sequence of $(1,\varepsilon_i)$-pancakes where the measure $\nu$ is a (weak)-limit of the sequence of probability measures obtained by normalising the volume of each pancake.
\begin{de}[Constructing Pancake]
Let $(I,\nu)$ be a spherical needle. We call a $(1,\delta)$-pancake $S$, a constructing pancake for $(I,\nu)$, if there exists a decreasing sequence of pancakes $...\subset S_i\subset S_{i-1}\subset\cdots\subset S_0$, where $S_0=S$ and $(I,\nu)$ is a limit of this sequence.
\end{de}
\begin{de}[Distance Between Spherical Needles]
For $\varepsilon>0$, we say that the distance between the spherical needles $(I_1,\nu_1)$ and $(I_2,\nu_2)$ is at most equal to $\varepsilon$ if there exists a constructing pancake $S_1$ (resp $S_2$) for $(I_1,\nu_1)$ (resp $(I_2,\nu_2)$) such that $S_1$ (resp $S_2$) is in the $\varepsilon$-neighborhood of $I_1$ (resp $I_2$) and the Hausdorff distance between $S_1$ and $S_2$ is at most equal to $\varepsilon$.
\end{de}

\subsection{A Fundamental Spherical Localisation Lemma}
The main result of this section is the next lemma which is known as the \emph{localisation lemma}. The Euclidean counterpart is proved in \cite{lova}. The reader can skip the proof as it is very similar to the proof presented in \cite{lova}. 

\begin{lem} \label{genlova}
Let $G_{i}$ for $i=1,2$  be two continuous functions on $\mathbb{S}^{n}$ such that
\begin{eqnarray*}
\int_{\mathbb{S}^{n}}G_{i}(u)d\mu(u)>0,
\end{eqnarray*}
then a convex partition of $\mathbb{S}^n$, $\Pi\in \mathcal{CP}^{\leq 1}$ by spherical needles exists such that for every $\sigma$ an element of $\Pi$, we have
\begin{equation} \label{eqn:baba}
\int_{\sigma}G_{i}(t)d\nu_{\sigma}(t)>0.
\end{equation}
$\nu_{\sigma}$ is a $\sin^{n-1}$-affine probability measure which is canonically defined on $\sigma$ from the partition. 
\end{lem}
%\begin{lem} \label{fondaa}
%Let $f_1$, $f_2$, $f_3$, and $f_4$ be four continuous non-negative functions on $\mathbb{S}^n$. Suppose that for every $u\in \mathbb{S}^n$ we have $f_1(u).f_2(u)\leq f_3(u).f_4(u)$. Let $\mu$ denote the normalised Riemannian measure of $\mathbb{S}^{n}$. If for every probability measure $\nu$ having a $\sin^{n-1}$-affine density function and supported on a geodesic segment $I$, we have
%\begin{eqnarray*}
%\left(\int_{I} f_1(t) d\nu(t)\right)\left(\int_{I} f_2(t) d\nu(t)\right)\leq \left(\int_{I} f_3(t) d\nu(t)\right)\left(\int_{I} f_4(t) d\nu(t)\right),
%\end{eqnarray*}
%where the geodesic segment $I$ is parametrised by its arc length $t$ and
%\begin{eqnarray*}
%f(t)d\nu(t)=f(t)g(t)dt,
%\end{eqnarray*}
%$g(t)$ is the density of the measure $\nu$, then we will have:
%\begin{equation} \label{eqn:gen}
%\left(\int_{\mathbb{S}^{n}} f_1(u) d\mu(u)\right)\left(\int_{\mathbb{S}^{n}}f_2(u) d\mu(u)\right)\leq \left(\int_{\mathbb{S}^{n}} f_3(u) d\mu(u)%%\right)\left(\int_{\mathbb{S}^{n}} f_4(u) d\mu(u)\right).
%\end{equation}
%\end{lem}

%This Lemma is a spherical version of the Localisation Lemma and its corollary which are both proved in \cite{lova} in the Euclidean case. I shall first prove the following:
%\begin{lem} \label{lova}
%Let $G_{i}$ for $i=1,2$  be two continuous functions on $\mathbb{S}^{n}$ such that
%\begin{eqnarray*}
%\int_{\mathbb{S}^{n}}G_{i}(u)d\mu(u)>0,
%\end{eqnarray*}
%then a $\sin^{n-1}$-affine probability measure $\nu$ supported on a geodesic segment $I$ exists such that
%\begin{eqnarray*}
%\int_{I}G_{i}(t)d\nu(t)>0.
%\end{eqnarray*}
%\end{lem}

\emph{Proof of Lemma \ref{genlova}}

First step is to prove the following claim,

\emph{claim}:

There exists a spherical needle $(I,\nu)$ such that equation (\ref{eqn:baba}) is satisfied.

\emph{Proof of the claim}:

We construct a decreasing sequence of convex subsets of $\mathbb{S}^n$ using the following procedure:
\begin{itemize}
\item Define the first step cutting map $F_1:\mathbb{S}^n \to \mathbb{R}^2$ by
\begin{equation} \label{eqn:borsuk}
F_1(x)=\left(\int_{x^{\vee}}G_{1}(u)d\mu(u), \int_{x^{\vee}}G_{2}(u)d\mu(u)\right)
\end{equation}
where $x^{\vee}$ denotes the (oriented) open hemi-sphere centered at the point $x$. Apply Borsuk-Ulam Theorem to $F_1$. Hence there exists a $x_1^{\vee}$ such that
\begin{eqnarray*}
\int_{x_1^{\vee}}G_{1}(u)d\mu(u)=\int_{-x_1^{\vee}}G_{1}(u)d\mu(u) \\
\int_{x_1^{\vee}}G_{2}(u)d\mu(u)=\int_{-x_1^{\vee}}G_{2}(u)d\mu(u).
\end{eqnarray*}

Choose the hemi-sphere, denoted by $x^{\vee}_1$. Set $S_1=x^{\vee}_1\cap \mathbb{S}^{n}$.
\item Define the $i$-th step cutting map by
\begin{eqnarray*}
F_i(x)=\left(\int_{S_{i-1}\cap x^{\vee}}G_{1}(u)d\mu(u), \int_{S_{i-1}\cap x^{\vee}}G_{2}(u)d\mu(u)\right).
\end{eqnarray*}
By applying the Borsuk-Ulam Theorem to $F_i$, we obtain two new hemi-spheres and we choose the one, denoted by $x^{\vee}_i$. Set $S_i=x^{\vee}_i\cap S_{i-1}$.
\end{itemize}

This procedure defines a decreasing sequence of convex subsets $S_i=x^{\vee}_i\cap S_{i-1}$ for every $i\in\mathbb{N}$. Set: 
 
\begin{eqnarray*}
S_{\pi}=\bigcap_{i=1}^{\infty}(S_{i})=\bigcap_{i=1}^{\infty}clos(S_{i}),
\end{eqnarray*}
where $clos(A)$ determines the topological closure of the subset $A$. We call the hemi-spheres obtained from the cutting maps, the \emph{cutting hemi-spheres}.
\begin{de}[Cutting Hemi-spheres]
A cutting hemi-sphere is a $\mathbb{S}^n_{+}$ which is a hemi-sphere used at some stage of the algorithmic procedure described above. The first cutting hemi-sphere will be the hemi-sphere used at the very first stage of the procedure to cut the sphere $\mathbb{S}^n$ into two parts.
\end{de}

A convexly-derived probability measure $\nu_{\pi}$ is defined on $S_{\pi}$.
Since  $\lim_{i\to \infty}S_i=S_{\pi}$ (this limit is with respect to Hausdorff topology) the definition of the convexly-derived measures can be applied to define the positive probability measure supported on $S_{\pi}$ by
\begin{eqnarray*}
\nu_{\pi}=\lim_{i\to\infty}\frac{\mu|S_i}{\mu(S_i)}.
\end{eqnarray*}
Hence, by the definition of $\nu_{\pi}$
\begin{eqnarray*}
\int_{S_{\pi}}G_{j}(x)d\nu_{\pi}(x)=\lim_{i\to\infty}\frac{\displaystyle\int_{S_{i}}G_{j}(x)d\mu(x)}{\mu(S_{i})}
\end{eqnarray*}
for $j=1,2$, and where the limit is taken with respect to the vague topology defined on the space of convexly-derived measures (see \cite{memwst}). Recall the following :

\begin{lem}
\label{radon}
{\em (See \cite{hirsh}).}
Let $\mu_i$ be a sequence of positive Radon measures on a locally-compact space $X$ which vaguely converges to a positive Radon measure $\mu$. Then, for every relatively compact subset $A \subset X$, such that $\mu(\partial A)=0$,
$$\lim_{i \to \infty}\mu_i (A)=\mu(A).$$
\end{lem}

By the definition of the cutting maps $F_i(x)$, for every $i\in \mathbb{N}$, $j=1,2$ we have
\begin{eqnarray*}
\int_{S_{i}}G_{j}(u)d\mu(u)>0.
\end{eqnarray*}
By applying Lemma \ref{radon}, we conclude that the convexly-derived probability measure defined on $S_{\pi}$ satisfies the assumption of the Lemma \ref{genlova}. The dimension of $S_{\pi}$ is $<n$. Indeed, if it is not the case, then $dim S_{\pi}=n$. Since there is a convexly-derived measure with positive density defined on $S_{\pi}$, and by the construction of the sequence $\{S_i\}$ for every open set $U$ we have
\begin{eqnarray*}
\nu_{\pi}(S_{\pi}\cap U)&=&\lim_{i\to\infty}\frac{\mu(S_i\cap U)}{\mu(S_i)}\\
                          &=&\lim_{i\to\infty}\frac{\mu(S_{\pi}\cap U)}{2^i\mu(S_i)}.
\end{eqnarray*}
By supposition on the dimension of $S_{\pi}$, the right-hand equality is equal to zero. This is a contradiction with the positive measure $\nu_{\pi}$ charging mass on $S_{\pi}\cap U$.

Thus, $dim S_{\pi} <n$. If $dim S_{\pi}=1$ then the claim is proved. Note that $dim S_{\pi}$ can not be equal to zero, since the cutting map in each step cuts the set $S_{i}$ . If $dim S_{\pi}=k>2$, we define a new procedure by replacing $S$ with $S_{\pi}\cap S$, replacing the normalised Riemannian measure by the measure $\nu_{\pi}$, and replacing the sphere $\mathbb{S}^n$ by the sphere $\mathbb{S}^k$ containing $S_{\pi}$. For this new procedure, we define new cutting maps in every step. Since $k> 2$, by using the Borsuk-Ulam Theorem we obtain hyperspheres ($\mathbb{S}^{k-1}$) halving the desired (convexly-derived) measures. The new procedure defines a new sequence of convex subsets and, by the same arguments given before, a convexly-derived measure defined on the interersection of this new sequence satisfying the assumption of Lemma \ref{genlova}. By the same argument, the dimension of the intersection of the decreasing sequence of convex sets is $<k$. If the dimension of the intersection of this new sequence is equal to $1$,  we are finished. If not, we repeat the above procedure until arriving to a $1$-dimensional set. This proves that a probability measure $\nu$ with a (non-negative) $\sin^{n-1}$-concave density function $f$, supported on a geodesic segment $I$ exists such that:
\begin{equation} \label{eqn:dead}
\int_{I}f(t)G_i(t)dt\geq 0.
\end{equation}
We determine $I$ to have minimal length. If $f$ is $\sin^{n-1}$-affine on $I$ then we are done. We suppose this is not the case. We choose a subinterval $J\subset I$, maximal in length, such that a $\sin^{n-1}$-concave function $f$ satisfying $(\ref{eqn:dead})$ exists such that $f$ additionally is $\sin^{n-1}$-affine on the subinterval $J$. The existence of $J$ and $f$ follows from a standard compactness argument. We can assume that the length of $I$ is $<\pi/2$. Consider the Euclidean cone over $I$. Let $a,b\in I$ be the end points of $I$ and take the \emph{Euclidean segment} $[a,b]$ in $\mathbb{R}^2$ (basically the straight line joining $a$ to $b$). By definition of $\sin^{n-1}$-concave functions, the function $f$ is the restricion of a one-homogeneous $x^{n-1}$-concave function $F$ on the circle (a $x^{n-1}$-concave function $F$ is a function such that $F^{1/(n-1)}$ is concave). Transporting the entire problem to $\mathbb{R}^{n+1}$, we begin with two homogeneous functions $\bar{G_i}$ on $\mathbb{R}^{n+1}$ such that
\begin{eqnarray*}
\int_{\mathbb{R}^{n+1}}\bar{G_i} dx >0
\end{eqnarray*}
and we proved that there is a $2$-dimensional cone over a segment $[a,b]$, a one-homogeneous $x^{n-1}$-concave function $F$ on $[a,b]$, and a subinterval $[\alpha,\beta]\subset [a,b]$ such that $F^{1/(n-1)}$ is linear on $[\alpha,\beta]$ (this is due to the fact that by definition, the restriction of a one-homogeneous $x^{n-1}$-affine function on a $2$-dimensional Euclidean cone defines a $\sin^{n-1}$-affine function on the circle) and such that
\begin{eqnarray*}
\int_{[a,b]}\bar{G_i}(t)F(t)dt\geq 0.
\end{eqnarray*}
We echo the arguments given in \cite{lova} (pages $21-23$) (with the only difference being that every construction there drops by one dimension). This drop of dimension is necessary so that that every construction may preserve homogeneity- or in other words, one dimension must be preserved for the $2$-dimensional cone defined on $[a,b]$. Hence the proof of the claim follows.
\begin{flushright}
$\Box$
\end{flushright}

In the proof of the claim we used a family of $\{x^{\vee}_i\}$ of oriented hemi-spheres to cut the sphere. Each $x^{\vee}_{i}$ cuts the sphere in two parts in such a way that in \emph{both} parts of the sphere the integral of $G_i$ remains positive (due to the Borsuk-Ulam Theorem). At each stage of the cutting, we only kept one part of the sphere. However, if we carry out everything we did with respect to the \emph{other} parts, we obtain (in a straightforward way) the conclusion of Lemma \ref{genlova}.

\begin{flushright}
$\Box$
\end{flushright}  

\emph{Remark}:
\begin{itemize}
\item We should keep in mind that the partition $\Pi\in \mathcal{CP}^{\leq 1}$ can be constructed by choosing the family of cutting hemi-spheres $\{x^{\vee}_i\}$ such that all the vectors $x_i$ belong to a sphere of dimension $k$ provided $k\geq 2$. This has two benefits:
\begin{itemize}
\item The partition obtained in Lemma \ref{genlova} is \emph{not} unique.
\item We can choose the \emph{direction} of the cuts by appropriately choosing the sphere $\mathbb{S}^k$. This fact will become very useful in the proof of Theorem \ref{strong}.
\end{itemize} 
\item Instead of applying the Borsuk-Ulam Theorem in the proof of Lemma \ref{genlova} we can use the more powerful Gromov-Borsuk-Ulam Theorem which is stated and proved in \cite{memwst}. The Gromov-Borsuk-Ulam Theorem provides a convex partition $\mathcal{CP}^{\leq k}$ for every continuous map $f:\mathbb{S}^n\to\mathbb{R}^k$ where $k<n$ and a point $z\in\mathbb{R}^k$ such that $f^{-1}(z)$ intersects the maximum points of the density of the convexly-derived-measures associated to the partition. Therefore by applying the Gromov-Borsuk-Ulam Theorem directly for the map $f:\mathbb{S}^n\to \mathbb{R}^2$ defined by
\begin{eqnarray*}
f(x)=\left(\int_{x^{\vee}}G_{1}(u)d\mu(u), \int_{x^{\vee}}G_{2}(u)d\mu(u)\right),
\end{eqnarray*}
we obtain the desired convex partition $\Pi\in\mathcal{CP}^{\leq 1}$ of the Lemma \ref{genlova}.
\end{itemize}

\section{Some Anisotropic $2$-Dimensional Correlation Problems}

In this section, some correlation problems in $\mathbb{R}^2$ will be presented. The ideas from this section are to be compared with those of \cite{figa} (pointed out by M.Ledoux). The notations used in this section will be the same as in \cite{figa}. The main goal of this long section is to prove Corollary \ref{corrr}. In order to achieve this goal some definitions, lemmas and theorems will be needed.

\begin{de}[Strips]
A set $S\subset \mathbb{R}^2$ is called a strip if $S$ is open, convex, symmetric with respect to the origin, and if a $u\in \mathbb{S}^1$ and $h>0$ exist such that 
\begin{eqnarray*}
S=\{x\in \mathbb{R}^2 : \vert x.u\vert<h\}.
\end{eqnarray*}
$h$ is the \emph{width} of the strip and $u$ is the unit vector of the axis of the strip. The angle $\theta$ of two strips $S$ and $S'$ is equal to the angle between the respective unit vectors of the axis of the strips.
\end{de}

\begin{de}[Angular-Length function]
Let $E$ be an open set containing the origin. The function $\theta_{E}:(0,\infty)\to [0,\pi/2]$ is defined as
\begin{eqnarray*}
\theta_{E}(r)=\frac{1}{4}\frac{H^1(E\cap \partial B_r)}{r},
\end{eqnarray*}
where $H^1$ stands for the $1$-dimensional Hausdorff measure.
\end{de}
It is clear that for every $r>0$, $\theta_{E}(r)\leq \pi/2$.

\begin{de}[Width Decreasing Sets]
A set $E\subset \mathbb{R}^2$ is said to be \emph{width-decreasing} if $E$ is open, contains the origin and is symmetric with respect to it, and for every $r>0$ if $\theta_{E}(r)<\pi/2$ then $\theta_{E}\leq \theta_{S}$ on $(r,\infty)$ where $S$ is any strip for which $\theta_{E}(r)=\theta_{S}(r)$.
\end{de}

\emph{Remark}:

In \cite{figa}, it is shown that every symmetric convex body in $\mathbb{R}^2$ is a width-decreasing set. This fact will be used throughout this section. However, every width-decreasing set is not necessarily a symmetric convex set.

The following definition will be useful in the proof of Corollary \ref{rot} and Corollary \ref{corrr}:
\begin{de}[Nice Convex Bodies] \label{nicec}
Let $K$ be a symmetric convex set in $\mathbb{R}^2$. We call the set $K$ a ``\emph{nice convex body}'' if for every $r>0$ such that $\partial B(0,r)\cap K\neq 0$, there exists an arc $I\subset \partial B(0,r)\cap K$ such that $pr(I)\subset \mathbb{S}^1$ is an arc which contains the point $(-1,0)$ (or the point $\theta=\pi$ on the circle). Here $pr$ is the radial projection map for every $r>0$, maps $\partial B(0,r)\cap K$ to $\mathbb{S}^1$. 
\end{de}

\emph{Example}: The symmetric caps and symmetric strips are two examples of \emph{nice} convex bodies.

In this section, we will be working in $\mathbb{R}^2$ and we will be using a family of anisotropic measures which will be denoted by $\mu_{2,\beta}$. Before beginning, the definition of these measures will be given, and in the next section we shall see the reason why it's necessary to deal with such  unusual measures.
 
\begin{de}[The measure $\mu_2$] \label{mes}

Let $f$ be a continuous radial function defined on $\mathbb{R}^2$. Denote $\mathbb{S}^1_{+}$ by the half unit circle. Take a $\sin^{n-1}$-affine measure of the form $d\nu(t)=g(t)dt=C(n)\cos(t)^{n-1}dt$. The support of $\nu$ is a half circle $\mathbb{S}^1_{+}$. Extend this measure on the whole $\mathbb{S}^1$ in such a way that for every $u\in \mathbb{S}^1$, $g(u)=g(-u)$ where $g$ is the density function of this measure. From now on, we are supposing that the $\sin^{n-1}$-concave measures and functions are extended to $\mathbb{S}^1$. The measure $\mu_2$ is the measure defined as $r^n f(r)dr\wedge d\nu(t)$ in polar coordinates of $\mathbb{R}^2$. In Cartesian coordinates $x\,y$, the measure $\mu_2$ is equal to the measure:
\begin{eqnarray*}
\mu_2=C(n)\vert y\vert^{n-1}f(\sqrt{(x^2+y^2)}dx dy,
\end{eqnarray*}
here, the end point of the unit vector of the $y$-axis coincides with the maximum point of the function $g$. The constant $C(n)$ is a normalisation constant which is defined to be:
\begin{eqnarray*}
C(n)=(\int_{-\pi/2}^{+\pi/2}\cos(t)^{n-1}dt)^{-1}.
\end{eqnarray*}
\end{de} 
\emph{Remark}: The normalisation constant $C(n)$ does not play an important role since it cancels out within the inequalities we are dealing with.
Note also that the measure $\mu_2$ is \emph{not} rotationnally invariant, and therefore it makes sense to define new measures from $\mu_2$ by considering a \emph{rotation of the spherical part} of this measure.

\begin{de}[The Measures $\mu_{2,\beta}$]
Let $0\leq \beta\leq \pi$. The measure $\mu_{2,\beta}$ is derived from the measure $\mu_2$ by rotating the density function of the measure of the spherical needle by the angle $\beta$. Therefore $\mu_{2,\beta}=C(n) r^n f(r)\cos(t-\beta)^{n-1}dr\wedge dt$. The cartesian coordinates associated to the measure $\mu_{2,\beta}$ are the cartesian coordinates in $\mathbb{R}^2$ where the end-point of the unit vector of the $y$-axis in these coordinates coincides with the maximum point of the density function $\cos(t-\beta)^{n-1}$. 
\end{de}  

\emph{Remark}: One could similarily define the measure by considering $\cos(t+\beta)^{n-1}$. However, since we shall always be dealing with symmetric convex bodies in $\mathbb{R}^2$ these measures would lead to the same result as for $\cos(t-\beta)^{n-1}$.

\subsection{Symmetrisation of Symmetric Convex Sets Intersecting the Half-Plane}

The \emph{cap} symmetrisation process is a useful process in convex geometry which enables one to simplify the study of some geometric problems, notably of isoperimetric type. This symmetrisation process will be used throughout the rest of this section. The necessary definitions will be given below.

%\emph{Remark}:

%Before beginning the proof of this theorem, let us examine where the measure $\mu_2$ is coming from. If for all $y\geq 0$, we write down this measure in polar coordinates, then:
%\begin{eqnarray*}
%\mu_2&=&C(n)y^{n-2}e^{-(x^2+y^2)/2}dx\wedge dy \\
%     &=&C(n) r^{n-1}e^{-r^2/2}\sin(\theta)^{n-2}dr\wedge d\theta,
%\end{eqnarray*}     
%where $0\leq \theta \leq \pi$. As we can see now, the measure $\mu_2$ is the measure $\mu_1$ where one replaces the $\sin^{n-2}$-concave function $g(\theta)$ with the \emph{model} function $\sin(\theta)^{n-2}$. Additionally, note that the maximum point of $g(\theta)$ and $\sin(\theta)^{n-2}$ coincides. We denote by $\nu_2=C(n)\sin(\theta)^{n-2}d\theta$.

Let the $xy$-coordinate be such that the point $(0,1)$ coincides with the maximum point of $g(t)$. We assume $\mathbb{S}^1$ is parametrised canonically by $t$. Therefore $t=0$ corresponds to $(1,0)$ and $t=\pi/2$ corresponds to the maximum point of the function $g$. For every $r>0$ let $I_{K_i}(r)=pr(\partial B(0,r)\cap K_i)$ for $i=1,2$ and where $pr$ is the radial projection of $\partial B(0,r)$ to $S^1$. Let $\mathcal{M}_0$ be the set of all open sets containing the origin and symmetric with respect to the origin. Let $\alpha\in \mathbb{S}^1$ and $E\in \mathcal{M}_0$.

\begin{de}[$\alpha$-Cap Symmetrisation]
The map $s_{\alpha}:\mathcal{M}_0\to\mathcal{M}_0$ which we call the $\alpha$-double cap symmetrisation map which is defined by
\begin{eqnarray*}
s_{\alpha}(E)=\cup_{r>0}(\{re^{i\phi} :\vert \phi-\alpha \vert\leq \varepsilon_{E}(r)\}\cup\{re^{i\phi} : \vert \phi-(\alpha+\pi)\vert\leq \varepsilon_{E}(r)\}),
\end{eqnarray*}
where $re^{i\phi}=(r\cos(\phi),r\sin(\phi))\in \mathbb{R}^2$ and $\varepsilon_{E}(r)$ is such that 
\begin{eqnarray*}
\int_{I_{E}(r)}g(t)dt=2\int_{\alpha-\varepsilon_{E}(r)}^{\alpha+\varepsilon_{E}(r)}g(t)dt.
\end{eqnarray*}
\end{de}

%\emph{Note}:
%Above we used the variable $\theta$ which corresponded to a point of $\mathbb{S}^1$. In what will follow, we may keep using this notation occasionally when we want to paramterise the circle $\mathbb{S}^1$. For the angular-length function of a set $E$, we shall use the notation $\theta_{E}$ which should not be confused with the variable $\theta$ explained above. And for the $\theta$-cap symmetrisation of a set $E$ with respect to a point $\theta$ on the circle, we will use the notation $s_{\theta}(E)$.

Unfortunately, the image of the set of convex sets (or width-decreasing sets) by $s_{\alpha}$ is not necessarily the set of convex sets (or width-decreasing sets). However one has the following:
\begin{lem} \label{cap}
Let $E$ be a symmetric convex set containing the origin. Let $\varepsilon$ be such that for every $r>0$ we have
\begin{eqnarray*}
\varepsilon_{E}(r)\leq \theta_{E}(r),
\end{eqnarray*}
where $\theta_{E}(r)$ is the angular-length function. Then $s_{\alpha}(E)$ is a width-decreasing set.
\end{lem}

\emph{Proof of Lemma \ref{cap}}

Define (as in \cite{figa}):
\begin{eqnarray*}
\varepsilon_{E}'(r)=\limsup_{\delta\to 0^{+}}\frac{\varepsilon_{E}(r+\delta)-\varepsilon_{E}(r)}{\delta},
\end{eqnarray*}
It is shown in \cite{figa} that in order for $s_{\alpha}(E)$ to be a width-decreasing set, it is sufficient to show that for every $r>0$ such that $\varepsilon_{E}(r)<\pi/2$, we have
\begin{eqnarray*}
\varepsilon_{E}'(r)\leq -\frac{\tan(\varepsilon_{E}(r))}{r}\leq 0.
\end{eqnarray*}

Let $r>0$ be such that $\theta_{E}(r)<\pi/2$ and let 
\begin{eqnarray*}
pr(\partial B(0,r)\cap E)=\cup_{i=1}^{N}I_i\cup J_i,
\end{eqnarray*}
where $I_i$ is a subarc of $\mathbb{S}^1$ and $J_i=\{-x\vert x\in I_i\}$.
\begin{eqnarray*}
\int_{I_{E}(r)}g(t)dt=\sum_{i=1}^{N}\int_{I_i}g(t)dt.
\end{eqnarray*}

For every $r>0$, we denote:
\begin{eqnarray*}
\partial B(0,r)\cap E=\cup_{i=1}^N T_i\cup (-T)_i,
\end{eqnarray*}
In \cite{figa}, it is shown that for every $1\leq i\leq N$ we have :
\begin{equation} \label{eqn:tst}
H^1(T_i^{\varepsilon})\leq H^1(T_i)+\frac{\varepsilon}{r}(H^1(T_i)-2r\tan(\theta_{E}(r)))+o(\varepsilon),
\end{equation}
as $\varepsilon\to 0$.
Hence :
\begin{eqnarray*}
\frac{H^1(T_i^{\varepsilon})}{r+\varepsilon}-\frac{H^1(T_i)}{r}&=&\frac{rH^{1}(T^{\varepsilon}_i)-rH^{1}(T_i)-\varepsilon H^1(T_i)}{r(r+\varepsilon)}\\
                                                             &\leq& \frac{-2\varepsilon\tan(\theta_{E}(r))}{r+\varepsilon}+\frac{o(\varepsilon)}{r+\varepsilon}.
\end{eqnarray*}
Therefore, as $\varepsilon\to 0$, using equation (\ref{eqn:tst}), we obtain that for every $r>0$ such that $\theta_{E}(r)<\pi/2$ and for every $1\leq i \leq N$ we have :
\begin{eqnarray*}
l(I_i)(r)'\leq \frac{-2\tan(\theta_{E}(r))}{r},
\end{eqnarray*}
where $l$ is the arc length in $\mathbb{S}^1$. In other words, if we assume the arc $I_i$ be paramterised canonically on $\mathbb{S}^1$ by the interval $[x-\tau_1,x+\tau_1]$, where $\tau_1$ is a function of $r$. Therefore we obtain: 
\begin{equation} \label{eqn:tau1}
\tau_1(r)'\leq \frac{-\tan(\theta_{E}(r))}{r}.
\end{equation}

Let $P=(0,1)\in\mathbb{R}^2$. Since the function $g(t)$ is assumed to be $\sin^n$-affine, without loss of generality we can suppose that the circle $\mathbb{S}^{1}-\{P\}$ is the interval $[-\pi,\pi]$ and the function $g(t)=\cos(t)^n$. The maximum point of the function $g(t)$ corresponds to the point $P=(0,1)\in\mathbb{R}^2$ which also corresponds to $t=0$ in the interval $[-\pi,\pi]$.

There exist $a,a_i,b_i$ for $i\in\{1,\cdots,N\}$ such that
\begin{eqnarray*}
\int_{a-\varepsilon(r)}^{a+\varepsilon(r)}\cos(t)^{n}dt=\sum_{i=1}^{N}\int_{a_i}^{b_i}\cos(t)^{n}dt.
\end{eqnarray*}
\begin{itemize}
\item For every $r>0$ if there exists an interval $I_j=[a_j,b_j]$ such that
\begin{equation} \label{eqn:three}
\vert\cos(a_j)^n-\cos(b_j)^n\vert\geq\vert\cos(\varepsilon-a)^n-\cos(\varepsilon+a)^n\vert,
\end{equation}
we get:
\begin{eqnarray*}
\varepsilon_{E}'(r)&\leq& \theta_{E}'(r)\\
               &\leq&-\frac{\tan(\theta_{E}(r))}{r}\\
               &\leq&-\frac{\tan(\varepsilon_{E}(r))}{r}.
\end{eqnarray*}

\item For every $r>0$ such that such an interval $I_j$ (for every $j\in\{1,\cdots,N\}$) satisfying the above inequality of ($\ref{eqn:three}$) does not exist, thanks to the $\sin^n$-concavity of $\cos(x)^n$ (and to $\theta\geq \varepsilon$), we obtain:
\begin{eqnarray*}
\varepsilon'(r)&\leq&\frac{\displaystyle\sum_{i=1}^{N}\theta_{E}'(r)\vert\cos(a_i)^n-\cos(b_i)^n\vert}{\vert\cos(a+\varepsilon)^n-\cos(a-\varepsilon)^n\vert}\\
               &\leq&-\frac{\displaystyle\sum_{i=1}^{N}N\tan(\theta_{E}(r))\vert\cos(a_i)^n-\cos(b_i)^n\vert}{r\vert\cos(a+\varepsilon)^n-\cos(a-\varepsilon)^n\vert}\\
               &\leq&-\frac{N\tan(\theta_{E}(r))}{r}\\
               &\leq&-\frac{\tan(\varepsilon_{E}(r))}{r}.
\end{eqnarray*}
\end{itemize}

This ends the proof of the Lemma \ref{cap}.
\begin{flushright}
$\Box$
\end{flushright}

\emph{Remark}:
\begin{itemize}
\item According to Lemma \ref{cap} and Lemma \ref{sinconc}, for the vertical double cap symmetrisation (i.e. $s_{\pi/2}$), we are certain that the image of a convex set under $s_{\pi/2}$ is always a width-decreasing set. 
\item We shall be using equation (\ref{eqn:tau1}) frequently whenever we have to check if any set has the width-decreasing property.
\end{itemize}
\subsection{An Algorithmic Procedure Assigned to Two Symmetric Convex Sets}

%Returning to the proof, we apply the vertical double cap symmetrisation map, $s_{\pi/2}$ to $K_i$ for $i=1,2$. Let $r_i$ be the radius of the largest disk inscribed in $K_i$. By definition, it is clear that for $r\geq r_i$, we have
%\begin{eqnarray*}
%\theta_{K_i}(r)\leq \theta_{S_i}(r),
%\end{eqnarray*}
%where $S_i$ is a strip of width equal to $r_i$.

Let $K_1, K_2$ be two symmetric convex bodies. Remember that the measure with which we work is not rotationally  invariant, so the measure of strips of equal width but different axis are all different. Since we adjust the maximum of the density to be on the point $(1,0)$, the vertical strip of width equal to $r_i$ has the largest measure amoung other strips of width equal to $r_i$. 

We define a family of width-decreasing sets $\{s_{\pi/2}(S_{\alpha,r_i})\}_{\alpha}$ parametrised by $\pi/2\leq \alpha \leq\pi$. This is precisely the family of vertical double cap symmetrisation of strips with width equal to $r_i$, where the unit vector of the axis of the strip varies from $\alpha=\pi$ to $\pi/2$. It is clear by Lemma \ref{sinconc} that for $r>0$ and $\pi/2\leq \alpha_1\leq \alpha_2\leq \pi$ we have $s_{\pi/2}(S_{\alpha_2,r})\subseteq s_{\pi/2}(S_{\alpha_1,r})$.

For $i=1,2$, let $r_i$ be the radius of the largest disk inscribed in $K_i$. Define now the following algorithmic procedure for $K_1$ and $K_2$:
\begin{itemize}
\item First step: If $s_{\pi}(K_i)$ is width decreasing for $i=1$ or $i=2$, complete this procedure by performing the following symmetrisation : $s_{\pi}(K_i)$ and $s_{\pi/2}(K_j)$ for $i\neq j$. Otherwise, go to the next step.
\item Second step: Set $\alpha_i$ to be the largest $\pi/2\leq \alpha\leq \pi$ such that $s_{\pi/2}(K_i)\subset s_{\pi/2}(S_{\alpha,r_i})$ for $i=1,2$. Perform the following symmetrisation : $s_{\alpha_i}(K_i)$ and $s_{\pi/2}(K_j)$ for $i\neq j$. If 
\begin{eqnarray*}
\mu_2(K_1\cap K_2)\geq \mu_2(s_{\pi/2}(K_1)\cap s_{\alpha_2}(K_2)),
\end{eqnarray*}
end this procedure. Otherwise, go to the next step.
\item Third step: Set $\alpha_j$ to be the smallest $0\leq \alpha\leq \pi/2$ such that $s_{\pi/2}(K_j)\subset s_{\pi/2}(S_{\alpha_j,r_j})$. Perform the following symmetrisation : $s_{\alpha_i}(K_i)$ and $s_{\alpha_j}(K_j)$ for $i\neq j$. Stop here.
\end{itemize}

We denote the results obtained after performing the above procedure by $s_i(K_i)$ for $i=1,2$.

It is clear according to Lemma \ref{cap}, that in every step the symmetric sets defined after the cap symmetrisation operation remain width-decreasing sets. What is important to us is to compare the measure of the intersection \emph{after} the symmetrisation procedure, and the measure of the intersection \emph{before} this procedure. This is provided by the following:

\begin{lem} \label{inter}
For every step of the previous procedure, we have:
\begin{eqnarray*}
\mu_2(K_1\cap K_2)\geq \mu_2(s_1(K_1)\cap s_2(K_2))
\end{eqnarray*}
\end{lem}

\emph{Proof of Lemma \ref{inter}}:

Suppose the procedure ends at the first step. In this case, the lemma is proved by the inclusion-exclusion principle. Indeed, in case for a $r>0$ we have
\begin{eqnarray*}
\nu_1(s_1(I_1(r))\cap s_2(I_2(r)))\neq 0,
\end{eqnarray*}
reminding that $\nu_1(I_i(r))=\nu_1(s_i(I_i(r))$, we have
\begin{eqnarray*}
\nu_1(s_1(I_1(r))\cap s_2(I_2(r)))&=&\nu_1(s_1(I_1(r))+\nu_1(s_2(I_2(r))-1 \\
                                  &\leq& \nu_1(I_1(r))+\nu_1(I_2(r))-\nu_1(I_1(r)\cup I_2(r))\\
                                  &=& \nu_1(I_1(r)\cap I_2(r)).
\end{eqnarray*}                                  
This proves the lemma for the first procedure. 

For the second procedure, the lemma is settled by definition. It remains to prove this lemma for the third step.

It is clear that if $\alpha_j=\pi/2-\alpha_i$ then the proof is similar to the first step via the inclusion-exclusion principle. Suppose then $\theta_i-\theta_j\leq \pi/2$ and suppose by contradiction that we have:
\begin{eqnarray*}
\mu_2(K_1\cap K_2)\leq \mu_2(s_1(K_1)\cap s_2(K_2)).
\end{eqnarray*}

This means that there exists a $r>0$ such that $\varepsilon_1(r)+\varepsilon_2(r)\leq \pi/2$ and such that 
\begin{eqnarray*}
\nu_1(I_1(r)\cap I_2(r))\leq \nu_1(s_1(I_1(r))\cap s_2(I_2(r))).
\end{eqnarray*}
Therefore two sets $J_r$ and $J'_r$ exist such that 
\begin{eqnarray*}
J_r\subset s_1(I_1(r))\cap s_2(I_2(r)),
\end{eqnarray*}
and 
\begin{eqnarray*}
J'_r\subset (s_1(I_1(r))\cup s_2(I_2(r)))^c,
\end{eqnarray*}
where $A(r)^c$ is the complementary of the set $A(r)$ in $\partial B(0,r)$ and $J'_r\subset I_1(r)$ and $J'_r\not\subset I_2(r)$ and $\nu_1(J_r)=\nu_1(J'_r)$. This fact shows that we can find $\alpha< \alpha_1$ such that the cap symmetrisation $s_{\alpha}(K_1)$ is a width-decreasing set, which is a contradiction with the definition of $\alpha_1$. 
 
Then the proof of Lemma \ref{inter} follows.
\begin{flushright}
$\Box$
\end{flushright}

We are now prepared to prove the main result of this section:

\begin{theo} \label{twodim}
Let $\mathbb{R}^2$ be enhanced with the measure $\mu_2$ as in definition \ref{mes}. Let $K_1$ and $K_2$ be two centrally-symmetric convex bodies in $\mathbb{R}^2$. Then two symmetric strips $S_1$ and $S_2$ in $\mathbb{R}^2$ exist such that 
\begin{eqnarray*}
\frac{\mu_2(K_{1}\cap K_{2})}{\mu_2(K_{1})\mu_2(K_{2})}\geq \frac{\mu_2(S_{1}\cap S_{2})}{\mu_2(S_{1})\mu_2(S_{2})}.
\end{eqnarray*}
%where $\mu_2$ is the measure $C(n)\vert y\vert^{n-2}e^{-(x^2+y^2)/2}dxdy$ with respect to the $x-y$-coordinates in $\mathbb{R}^2$ such that the end point of the unit vector of the $y$-axis coincides with the maximum point of the function $g$. The constant $C(n)$ is defined to be:
%\begin{eqnarray*}
%C(n)=(\int_{-\pi/2}^{\pi/2}\cos(t)^{n-2}dt)^{-1}.
%\end{eqnarray*}
%$S_1$ and $S_2$ are either orthogonal, or the angle between them is smaller than $\pi/2$. 
\end{theo}

%\emph{Remark}:

%Before beginning the proof of this theorem, let us examine where the measure $\mu_2$ is coming from. If for all $y\geq 0$, we write down this measure in polar coordinates, then:
%\begin{eqnarray*}
%\mu_2&=&C(n)y^{n-2}e^{-(x^2+y^2)/2}dx\wedge dy \\
%     &=&C(n) r^{n-1}e^{-r^2/2}\sin(\theta)^{n-2}dr\wedge d\theta,
%\end{eqnarray*}     
%where $0\leq \theta \leq \pi$. As we can see now, the measure $\mu_2$ is the measure $\mu_1$ where one replaces the $\sin^{n-2}$-concave function $g(\theta)$ with the \emph{model} function $\sin(\theta)^{n-2}$. Additionally, note that the maximum point of $g(\theta)$ and $\sin(\theta)^{n-2}$ coincides. We denote by $\nu_2=C(n)\sin(\theta)^{n-2}d\theta$.

\emph{Proof of Theorem \ref{twodim}}

Apply the algorithmic procedure defined above to $K_i$ and denote the output by $s_i(K_i)$ for $i=1,2$.

Lemma \ref{inter} implies the following:

\begin{eqnarray*}
\frac{\mu_2(K_{1}\cap K_{2})}{\mu_2(K_{1})\mu_2(K_{2})}\geq \frac{\mu_2(s_1(K_1)\cap s_2(K_2)}{\mu_2(s_1(K_1))\mu_2(s_2(K_2))}.
\end{eqnarray*}

The rest of the proof is very similar to the one given in \cite{figa} but in the interest of being thorough, I will now break down the details.

Set
\begin{eqnarray*}
r_0= inf\{r>0 : \theta_{K_1}(r)+\theta_{K_2}(r)\leq \pi/2\},
\end{eqnarray*}
and without loss of generality, we could consider that $\theta_{K_1}(r_0)+\theta_{K_2}(r_0)\leq \pi/2$, $0<\theta_{K_1}(r_0)<\pi/2$ and $0<\theta_{K_2}(r_0)<\pi/2$.

Let $S_1$ be a cap such that the unit vector of its axis is given by $l_1$ and such that $\theta_{S_1}(r_0)=\theta_{K_1}(r_0)$. It is clear that
\begin{eqnarray*}
s_1(K_1)/B(0,r_0)\subset S_1\\
S_1\cap B(0,r_0)\subset s_1(K_1)\cap B(0,r_0).
\end{eqnarray*}
Define $S_2$ to be a strip such that the unit vector of its axis be given by $l_2$ and such that $\theta_{S_2}(r_0)=\theta_{K_2}(r_0)$. Similarily we have
\begin{eqnarray*}
s_2(K_2)/B(0,r_0)\subset S_2\\
S_2\cap B(0,r_0)\subset s_2(K_2)\cap B(0,r_0).
\end{eqnarray*}
Observe that if we set
\begin{eqnarray*}
E=(s_1(K_1)\cap B(0,r_0))\cup (S_1/B(0,r_0)) \\
F=(s_2(K_2)\cap B(0,r_0))\cup (S_2/B(0,r_0)).
\end{eqnarray*}
Then we have $E\cap F =s_1(K_1)\cap s_2(K_2)$ and clearly we obtain the following inequality
\begin{eqnarray*}
\frac{\mu_2(s_1(K_1)\cap s_2(K_2))}{\mu_2(s_1(K_1))\mu_2(s_2(K_2))}\geq \frac{\mu_2(E\cap F)}{\mu_2(E)\mu_2(F)}.
\end{eqnarray*}
And then
\begin{eqnarray*}
\frac{\mu_2(E\cap F)}{\mu_2(E)\mu_2(F)} &\geq& \frac{\mu_2(E\cap F)-\mu_2(E/S_1)}{(\mu_2(E)-\mu_2(E/S_1))(\mu_2(F))}\\
                                      &=& \frac{\mu_2(F\cap S_1)}{\mu_2(S_1)\mu_2(F)} \\
                                      &\geq& \frac{\mu_2(F\cap S_1)-\mu_2(F/S_2)}{\mu_2(S_1)(\mu_2(F)-\mu_2(F/S_2))}\\
                                      &=&\frac{\mu_2(S_1\cap S_2)}{\mu_2(S_1)\mu_2(S_2)}.
\end{eqnarray*}

This ends the proof of Theorem \ref{twodim}.
\begin{flushright}
$\Box$
\end{flushright}

\emph{Remark}:

Even if the result of Theorem \ref{twodim} is interesting in the sense that it provides us with two symmetric strips, there is no control on the axis of these symmetric strips. It turns out (and we shall see this later on in the next section) that having a control on the axis of at least one of the symmetric strips is in fact crucial for us. The next lemma asserts that for at least one $\beta$, with respect to the measure $\mu_{2,\beta}$, we can be certain of the position of the axis of at least one of these symmetric strips.

\begin{cor} \label{rot}
Let $\mathbb{R}^2$ be enhanced with the measure $\mu_2$, as defined in Theorem \ref{twodim}. Let $K_1$ and $K_2$ be two centrally-symmetric convex bodies in $\mathbb{R}^2$. There exists a measure $\mu_{2,\beta}$ obtained from the measure $\mu_2$ such that :
\begin{eqnarray*}
\frac{\mu_{2,\beta}(K_1\cap K_2)}{\mu_{2,\beta}(K_1)\mu_{2,\beta}(K_2)}\geq \frac{\mu_{2,\beta}(S_1\cap S_2)}{\mu_{2,\beta}(S_1)\mu_{2,\beta}(S_2)},
\end{eqnarray*}
where $S_1$ and $S_2$ are symmetric strips and the axis of either $S_1$ or $S_2$ coincides with the $x(\beta)$-axis, where $x(\beta)$ is the Cartesian coordinates assigned to the measure $\mu_{2,\beta}$ which is obtained by the rotation of the axis $x$ by the angle $\beta$.
\end{cor}
\emph{Remark}:

This corollary is a consequence of the continuity of the density functions of the family of measures $\mu_{2,\beta}$.

\emph{Proof of Corollary \ref{rot}} :

Define first:
\begin{eqnarray*}
Max(K_1)=\max_{r\geq 0}(\vert J_r(K_1)\vert)
\end{eqnarray*}
where $J_r(K_1)=\partial B(0,r)\cap K_1$ and $\vert J_r(K_1)\vert$ is the number of arcs obtained by intersecting with $\partial B(0,r)$. And define:
\begin{eqnarray*}
Max(K_1,K_2)=Max\{Max(K_1),Max(K_2)\}.
\end{eqnarray*}

%We prove Lemma \ref{rot} by induction on $Max(K_1,K_2)$. 

If $Max(K_1,K_2)=2$ (which is the possible minimum obtained due to the symmetry of $K_1$ and $K_2$), the proof of the Lemma follows by simply rotating the $x$-axis to the $x(\beta)$-axis for which either $K_1$ or $K_2$ will be a \emph{nice convex body}, and then by performing the symmetrisation used in the proof of Theorem \ref{twodim} to obtain the appropriate symmetric strips $S_1$ and $S_2$.

%We suppose that we can prove the Lemma for every $K_1$ and $K_2$ with $Max(K_1,K_2)=2N>2$. Suppose now we have $K_1$ and $K_2$ with $Max(K_1,K_2)=2N+2$. Without loss of generality we can assume that $Max(K_1)=2N+2$. Let $r_{cric}>0$ be the minimum of $r>0$ such that $\vert J_r(K_1)\vert=2N+2$.

Suppose then $Max(K_1,K_2)>2$. Following the arguments of Theorem \ref{twodim}, for every $\beta$, there exist two symmetric strips $M_{1,\beta}$, $M_{2,\beta}$ such that $Max(M_{1,\beta},M_{2,\beta})=2$ and such that
\begin{equation} \label{eqn:hihi}
\frac{\mu_{2,\beta}(K_1\cap K_2)}{\mu_{2,\beta}(K_1)\mu_{2,\beta}(K_2)}\geq \frac{\mu_{2,\beta}(M_{1,\beta}\cap M_{2,\beta})}{\mu_{2,\beta}(M_{1,\beta})\mu_{2,\beta}(M_{2,\beta})}.
\end{equation}

Of course, when $\beta$ changes in the equation (\ref{eqn:hihi}) , the sets $M_{1,\beta}$ and $M_{2,\beta}$ may also be moving. 

For $\beta=0$, let $M_1$ and $M_{2,0}$ be two width-decreasing sets such that $Max(M_1,M_{2,0})=2$ and such that :
\begin{eqnarray*}
\frac{\mu_{2,0}(K_1\cap K_2)}{\mu_{2,0}(K_1)\mu_{2,0}(K_2)}\geq \frac{\mu_{2,0}(M_{1}\cap M_{2,0})}{\mu_{2,0}(M_{1})\mu_{2,0}(M_{2,0})}. 
\end{eqnarray*}
This of course can be given by applying Theorem \ref{twodim}.

For $\beta>0$, define the constant $c(\beta)>0$ such that :
\begin{eqnarray*}
\mu_{2,\beta}(K_1)=c(\beta)\mu_{2,\beta}(M_1),
\end{eqnarray*}
and use the $\alpha$-cap symmetrisation to obtain a width-decreasing set $M_{2,\beta}$ such that $Max(M_{2,\beta})=2$, satisfying :
\begin{eqnarray*}
\mu_{2,\beta}(K_2)=\mu_{2,\beta}(M_{2,\beta}),
\end{eqnarray*}

Applying the arguments of Lemma \ref{inter}, we have:

\begin{eqnarray*}
\mu_{2,\beta}(K_1\cap K_2)\geq c(\beta)\mu_{2,\beta}(M_1\cap M_{2,\beta}).
\end{eqnarray*}
By continuity, the set $M_{2,\beta}$ can be constructed by (slightly) moving the set $M_{2,0}$. 

Hence :
\begin{equation} \label{eqn:thihi}
\frac{\mu_{2,\beta}(K_1\cap K_2)}{\mu_{2,\beta}(K_1)\mu_{2,\beta}(K_2)}\geq \frac{c(\beta)\mu_{2,\beta}(M_1\cap M_{2,\beta})}{c(\beta)\mu_{2,\beta}(M_1)\mu_{2,\beta}(M_{2,\beta})}.
\end{equation}

According to equation (\ref{eqn:thihi}), there exists a $\beta_0$ such that (with respect to the cartesian coordinates associated to the measure $\mu_{2,\beta_0}$) the set $M_1$ is a \emph{nice convex set}.

Therefore, applying Theorem \ref{twodim} for this $\beta_0$, the proof of the Corollary \ref{rot} follows.

\begin{flushright}
$\Box$
\end{flushright}

\subsection{Symmetrisation of Symmetric Convex Sets Intersecting General Cones}

Theorem \ref{twodim} and Corollary \ref{rot} alone are not sufficient for what will follow . Indeed, if we want to simplify an $n$-dimensional problem to a two-dimensional one by applying the spherical localisation, we end up with \emph{spherical needles} which can very well be segments of length strictly smaller than $\pi$. Therefore we are obliged to consider the intersection of our sets with different \emph{cones} in $\mathbb{R}^2$.

Suppose a spherical needle $I\subset \mathbb{S}^1_{+}$ be given. The cone defined by $I$ which is denoted by $C(I)$ is simply the cone in $\mathbb{R}^2$ with the vertex being the origin of $\mathbb{R}^2$, which contains the arc $I$, and the boundary of $C(I)$ contains the end points of the arc $I$. All the cones will be assumed to be provided by an arc $I$. For simplicity, a cone will be denoted by $C$.

For every $C$, we can consider the cone $-C$ and since we are dealing with symmetric convex bodies such as $K$, we have $C\cap K=-(-C\cap K)$

Let $K$ be a symmetric set in $\mathbb{R}^2$. Note that the angular length function of the set $(C\cup(-C))\cap K$ is defined similarily to the angular length function of $\mathbb{R}^2\cap K$. 

We say the set $(C\cup (-C))\cap K$ has the \emph{width-decreasing property} if the angular length function of this set satisfies the inequality :
\begin{eqnarray*}
\theta_{(C\cup(-C))\cap K}(r)'\leq \frac{-\tan(\theta_{(C\cup (-C))\cap K})}{r}.
\end{eqnarray*}

For simplifying the notations, one will omit working with $C\cup (-C)$ and only considers $C$. Indeed every construction one makes will be symmetric with respect to the origin and hence will work the same way for $-C$.

\begin{lem} \label{juju}
For every cone $C$ and symmetric convex body $K$, the set $K\cap C$ satisfies the \emph{width-decreasing propery}.
\end{lem}

\emph{Proof of Lemma \ref{juju}}:

For $r>0$, three different configurations will occur:
\begin{itemize}
\item For a certain $r>0$ we have $\partial B(0,r)\cap K=\partial B(0,r)\cap C$.
\item For a certain $r>0$ we have $\partial B(0,r)\cap K\neq \partial B(0,r)\cap C$ and $\partial B(0,r)\cap K$ contains (at least) one point on $\partial C$.
\item For a certain $r>0$ we have $\partial B(0,r)\cap K$ being strictly contained in the interior of $C$.
\end{itemize}

Clearly we only have to consider the second and third cases. If for a certain $r>0$, $\partial B(0,r)\cap K$ contains a point on the boundary of the cone, it means there exist at most two constant $a,b$ such that for every $r>0$ we have :
\begin{eqnarray*}
\theta_{K\cap \partial B(0,r)}=\cup_{i=1}^{2}[a,a_i]\cup J_r ,
\end{eqnarray*}
where $J(r)$ is every arc in $K\cap \partial B(0,r)$ inside the cone $C$ which does not intersect the boundary of $C$. Hence, by differentiating the above equality we get:
\begin{eqnarray*}
\theta'_{K\cap \partial B(0,r)}&=& \cup_{i=1}^{2}[a,a_i']\cup J'_r \\
                                 &\leq& \frac{-\tan(\theta_{K})}{r} \\
                                 &\leq& \frac{-\tan(\theta_{K\cap \partial B(0,r)})}{r}.
\end{eqnarray*}
(One used the inequality (\ref{eqn:tau1})). 

If for a certain $r>0$ we are in the configuration of the third case, then the fact that 
\begin{eqnarray*}
\theta'_{K\cap \partial B(0,r)}\leq \frac{-\tan(\theta_{K\cap \partial B(0,r)})}{r},
\end{eqnarray*}
is automatically verified since $K\cap \partial B(0,r)\subset K$ and $K$ is indeed a width-decreasing set.

This ends the proof of Lemma \ref{juju}.

\begin{flushright}
$\Box$
\end{flushright}

We now need to define the $\alpha$-cap symmetrisation process for the sets which are the intersection of a (symmetric) convex body with a cone. This will be slightly different from the one used in the proof of Theorem \ref{twodim}.

Let $\alpha\in[0,\pi]$. We define the $\alpha$-cap symmetrisation of $L=K\cap C$ as follows:
\begin{eqnarray*}
s_{\alpha}(L)=\cup_{r>0}(\{re^{i\phi} :\vert \phi-\alpha \vert\leq \varepsilon_{L}(r)\}\cup\{re^{i\phi} : \vert \phi-(\alpha+\pi)\vert\leq \varepsilon_{L}(r)\}),
\end{eqnarray*}
where $re^{i\phi}=(r\cos(\phi),r\sin(\phi))\in \mathbb{R}^2$ and $\varepsilon_{L}(r)$ is such that 
\begin{equation} \label{eqn:cone}
\frac{\displaystyle\int_{I_{L}(r)}g(t)dt}{\displaystyle\int_{\partial B(0,r)\cap C}g(t) dt}=\frac{\displaystyle\int_{\alpha-\varepsilon_{L}(r)}^{\alpha+\varepsilon_{L}(r)}g(t)dt}{\displaystyle\int_{-\pi/2}^{+\pi/2} g(t) dt},
\end{equation}
where $I_{L}(r)=pr(\partial B(0,r)\cap K\cap C)$ and $pr$ is the radial projection of $\partial B(0,r)$ to $\mathbb{S}^1$.

We are ready for the next theorem which is a conic version of theorem \ref{twodim}:

\begin{theo} \label{con}
Let $\mathbb{R}^2$ be enhanced with the measure $\mu_2$, as in definition \ref{mes}. Let $C$ be a cone in $\mathbb{R}^2$. Let $K_1$ and $K_2$ be two symmetric convex sets in $\mathbb{R}^2$. Then two strips $S_1$ and $S_2$ in $\mathbb{R}^2$ exist such that 
\begin{eqnarray*}
\frac{\mu_2(\mathbb{R}^2\cap C)\mu_2(K_1\cap K_2\cap C)}{\mu_2(K_1\cap C)\mu_2(K_2\cap C)}\geq \frac{\mu_2(\mathbb{R}^2)\mu_2(S_{1}\cap S_{2})}{\mu_2(S_{1})\mu_2(S_{2})}.
\end{eqnarray*}
\end{theo}

\emph{Proof of Theorem \ref{con}}:

If $C$ is the half-plane, the proof follows from Theorem \ref{twodim}. We (obviously) assume this is not the case, and $C$ is strictly contained in a half-plane. The trick here is to construct for each $K_i$, a width-decreasing set $M_i$ such that:
\begin{eqnarray*}
\frac{\mu_2(\mathbb{R}^2\cap C)\mu_2(K_1\cap K_2\cap C)}{\mu_2(K_1\cap C)\mu_2(K_2\cap C)}\geq \frac{\mu_2(\mathbb{R}^2)\mu_2(M_1\cap M_2)}{\mu_2(M_1)\mu_2(M_2)}.
\end{eqnarray*}

Let $L_i=C\cap K_i$. We parametrise the cone $C$. For $\rho_0<\pi/2$, let $C=C(0,2\rho_0)$ be a cone in $\mathbb{R}^2$ which is a cone over the arc $[x-\rho_0,x+\rho_0]\subset \mathbb{S}^1_{+}$ for a certain $x\in\mathbb{S}^1_{+}$. According to Theorem \ref{twodim}, there exists two symmetric strips $X_1$ and $X_2$ such that:
\begin{eqnarray*}
\frac{\mu_2(K_1\cap K_2)}{\mu_2(K_1)\mu_2(K_2)}\geq \frac{\mu_2(X_1\cap X_2)}{\mu_2(X_1)\mu_2(X_2)}.
\end{eqnarray*}
Denote the end-point of the unit vector of the axis of $X_i$ by $y_i$ (for $i=1,2$). Therefore $y_i\in\mathbb{S}^1$.

We perform the $\alpha$-cap symmetrisation for the set $C\cap K_1$ with respect to $y_1$.

\emph{claim}:
The set $s_{y_1}(C\cap K_1)$ is a width-decreasing set.

\emph{Proof of the claim}:

Let $\delta(r)$ be chosen such that for every $r>0$ we have:
\begin{eqnarray*}
\frac{\displaystyle\int_{I_{L_1}(r)}\cos(t)^{n-1} dt}{\displaystyle\int_{x-\rho_0}^{x+\rho_0}\cos(t)^{n-1} dt} =\frac{\displaystyle\int_{y_1-\delta}^{y_1+\delta}\cos(t)^{n-1} dt}{\displaystyle\int_{-\pi/2}^{+\pi/2}\cos(t)^{n-1} dt}.
\end{eqnarray*}
%Then we perform a symmetrisation of the cone $C$ and the set $L_1$ with respect to the new cone $C'=C'(0,2\rho_0)$ which is the cone over the segment $[-\rho_0,+\rho_0]\subset \mathbb{S}^1_{+}$ :

Let $\varepsilon(r)$ be such that:
\begin{eqnarray*}
\frac{\displaystyle\int_{-\varepsilon}^{+\varepsilon}\cos(t)^{n-1}dt}{\displaystyle\int_{-\rho_0}^{+\rho_0}\cos(t)^{n-1}dt}=\frac{\displaystyle\int_{I_{L_1}(r)}\cos(t)^{n-1} dt}{\displaystyle\int_{x-\rho_0}^{x+\rho_0}\cos(t)^{n-1} dt}.
\end{eqnarray*}

Therefore we have:
\begin{equation} \label{eqn:jaleb}
\frac{\displaystyle\int_{-\varepsilon}^{+\varepsilon}\cos(t)^{n-1}dt}{\displaystyle\int_{-\rho_0}^{+\rho_0}\cos(t)^{n-1}dt}=\frac{\displaystyle\int_{y_1-\delta}^{y_1+\delta}\cos(t)^{n-1} dt}{\displaystyle\int_{-\pi/2}^{+\pi/2}\cos(t)^{n-1} dt}.
\end{equation}

From equation (\ref{eqn:jaleb}) we can deduce that for every $r>0$, we have :
\begin{equation} \label{eqn:th}
\varepsilon(r)\leq \delta(r).
\end{equation}

By differentiating (\ref{eqn:jaleb}) we obtain:
\begin{eqnarray*}
2\varepsilon(r)'\cos(\varepsilon)^{n-1}=C\delta(r)'(\cos(y_1+\delta)^{n-1}+\cos(y_1-\delta)^{n-1}),
\end{eqnarray*}
where
\begin{eqnarray*}
C=\frac{\displaystyle\int_{-\rho_0}^{+\rho_0}\cos(t)^{n-1} dt}{\displaystyle\int_{-\pi/2}^{+\pi/2}\cos(t)^{n-1} dt}.
\end{eqnarray*}

According to Lemma \ref{juju}, the set $L_1$ satisfies the width-decreasing property hence:
\begin{eqnarray*}
\varepsilon(r)'\leq \frac{-\tan(\varepsilon)}{r}.
\end{eqnarray*}

Thus:
\begin{eqnarray*}
\delta(r)'&=&\frac{2\varepsilon(r)'\cos(\varepsilon)^{n-1}}{C(\cos(y_1+\delta)^{n-1}+\cos(y_1-\delta)^{n-1})}\\
          &\leq&\frac{-2\tan(\varepsilon)\cos(\varepsilon)^{n-1}}{rC(\cos(y_1+\delta)^{n-1}+\cos(y_1-\delta)^{n-1})}\\
          &\leq&\frac{-\tan(\varepsilon)\cos(\varepsilon)^{n-1}}{rC\cos(\delta)^{n-1}}.
\end{eqnarray*}
\begin{lem} \label{part}
The function 
\begin{eqnarray*}
\frac{\sin(x)\cos(x)^{n-2}}{\displaystyle\int_{0}^{x}\cos(t)^{n-1}dt},
\end{eqnarray*}
is decreasing on $[0,\pi/2]$.
\end{lem}

\emph{Proof of Lemma \ref{part}}:

Integrating by part, we obtain the following :
\begin{eqnarray*}
\displaystyle\int_{0}^{x}\cos(t)^{n-1}dt=\sin(x)\cos(x)^{n-2}+(n-2)\displaystyle\int_{0}^{x}\sin(t)^2\cos(t)^{n-3}dt.
\end{eqnarray*}
Since the function :
\begin{eqnarray*}
g(x)=\frac{\sin(x)^2\cos(x)^{n-3}}{\cos(x)^{n-1}}=\tan(x)^2,
\end{eqnarray*}
is an increasing function on $[0,\pi/2]$, therefore the function:
\begin{eqnarray*}
\frac{\displaystyle\int_{0}^{x}\sin(t)^2\cos(t)^{n-3}dt}{\displaystyle\int_{0}^{x}\cos(t)^{n-1}dt},
\end{eqnarray*}
is also an increasing function for $x\in[0,\pi/2]$ (here we used a well known Lemma which will be stated in the last section, Lemma \ref{bg}).
Therefore, we obtain :
\begin{eqnarray*}
1=\frac{\sin(x)\cos(x)^{n-2}}{\displaystyle\int_{0}^{x}\cos(t)^{n-1}dt}+\frac{(n-2)\displaystyle\int_{0}^{x}\sin(t)^2\cos(t)^{n-3}dt}{\displaystyle\int_{0}^{x}\cos(t)^{n-1}dt},
\end{eqnarray*}
and therefore we conclude that the function :
\begin{eqnarray*}
\frac{\sin(x)\cos(x)^{n-2}}{\displaystyle\int_{0}^{x}\cos(t)^{n-1}dt},
\end{eqnarray*}
is decreasing on $[0,\pi/2]$.
\begin{flushright}
$\Box$
\end{flushright}

Therefore according to Lemma \ref{part} and equation (\ref{eqn:th}), we obtain :
\begin{eqnarray*}
\frac{\sin(\varepsilon)\cos(\varepsilon)^{n-2}}{C} &=& \sin(\varepsilon)\cos(\varepsilon)^{n-2}\frac{\displaystyle\int_{-\delta}^{+\delta}\cos(t)^{n-1} dt}{\displaystyle\int_{-\varepsilon}^{+\varepsilon}\cos(t)^{n-1} dt} \\
&\geq& \sin(\delta).\cos(\delta)^{n-2},
\end{eqnarray*}
which in turn, means that:
\begin{eqnarray*}
\delta'(r)\leq \frac{-\tan(\delta)}{r}.
\end{eqnarray*}
Therefore the set $s_{y_1}(K_1\cap C)$ is a width-decreasing set in $\mathbb{R}^2$. This ends the proof of the claim.
\begin{flushright}
$\Box$
\end{flushright}

We perform the above $\alpha$-cap symmetrisation for the set $K_2\cap C$ with respect to $y_2\in\mathbb{S}^1$. Thus we obtain the following:
\begin{eqnarray*}
\frac{\mu_2(K_i\cap C)}{\mu_2(C)}=\frac{\mu_2(s_{y_i}(L_i))}{\mu_2(\mathbb{R}^2)},
\end{eqnarray*}
for $i=1,2$.
According to Lemma \ref{inter}, we have :
\begin{equation}\label{eqn:took}
\mu_2(K_1\cap K_2)\geq \mu_2(s_{y_1}(K_1)\cap s_{y_2}(K_2)),
\end{equation}
and therefore :
\begin{equation} \label{eqn:pen}
\frac{\mu_2(K_1\cap K_2\cap C)}{\mu_2(C)}\geq \frac{\mu_2(s_{y_1}(L_1)\cap s_{y_2}(L_2))}{\mu_2(\mathbb{R}^2)}.
\end{equation}
Combining equations (\ref{eqn:took}) and (\ref{eqn:pen}), we obtain :
\begin{eqnarray*}
\frac{\mu_2(C)\mu_2(K_1\cap K_2\cap C)}{\mu_2(K_1\cap C)\mu_2(K_2\cap C)}\geq \frac{\mu_2(\mathbb{R}^2)\mu_2(s_{y_1}(L_1)\cap s_{y_2}(L_2))}{\mu_2(s_{y_1}(L_1))\mu_2(s_{y_2}(L_2))}.
\end{eqnarray*}

The proof of Theorem \ref{con} then follows from Theorem \ref{twodim}.

\begin{flushright}
$\Box$
\end{flushright}

The following corollary is the counter-part of Corollary \ref{rot} when we intersect with general cones. In fact, the next corollary is the most general result which we shall need for the proof of Theorem \ref{strong}.

\begin{cor} \label{corrr}
Let $\mathbb{R}^2$ be enhanced with the measure $\mu_2$ as in definition \ref{mes}. Let $C(\beta)$ be a family of cones in $\mathbb{R}^2$ parametrised by $\beta \in [0,\pi]$. Let $K_1$ and $K_2$ be two symmetric convex sets in $\mathbb{R}^2$. Then a $\beta_0\in[0,\pi]$ and two symmetric strips $S_1$ and $S_2$ in $\mathbb{R}^2$ exist such that 
\begin{equation} \label{eqn:shahkar}
\frac{\mu_{2,\beta_0}(\mathbb{R}^2\cap C(\beta_0))\mu_{2,\beta_0}(K_1\cap K_2\cap C(\beta_0))}{\mu_{2,\beta_0}(K_1\cap C(\beta_0))\mu_{2,\beta_0}(K_2\cap C(\beta_0))}\geq \frac{\mu_{2,\beta_0}(\mathbb{R}^2)\mu_{2,\beta_0}(S_{1}\cap S_{2})}{\mu_{2,\beta_0}(S_{1})\mu_2(S_{2})}.
\end{equation}
Furtheremore, the axis of either $S_1$ or $S_2$ coincides with the $x(\beta_0)$-axis of the measure $\mu_{2,\beta_0}$.
\end{cor}

\emph{Proof of Corollary \ref{corrr}} :

The proof of Corollary \ref{corrr} is similar to the proof of Lemma \ref{rot}.

Following the arguments of Theorem \ref{con}, for every $\beta$, there exist two symmetric strips $M_{1,\beta}$, $M_{2,\beta}$ such that $Max(M_{1,\beta},M_{2,\beta})=2$ and such that
\begin{equation} \label{eqn:lhihi}
\frac{\mu_{2,\beta}(C(\beta)\cap \mathbb{R}^2)\mu_{2,\beta}(C(\beta)\cap K_1\cap K_2)}{\mu_{2,\beta}(C(\beta)\cap K_1)\mu_{2,\beta}(C(\beta)\cap K_2)}\geq \frac{\mu_{2,\beta}(\mathbb{R}^2)\mu_{2,\beta}(M_{1,\beta}\cap M_{2,\beta})}{\mu_{2,\beta}(M_{1,\beta})\mu_{2,\beta}(M_{2,\beta})}.
\end{equation}

Similar to Lemma \ref{rot}, when $\beta$ changes in equation (\ref{eqn:thihi}), the sets $M_{1,\beta}$ and $M_{2,\beta}$ may be moving. 

For $\beta=0$, let $M_1$ and $M_{2,0}$ be two width-decreasing sets such that $Max(M_1,M_{2,0})=2$ and such that :
\begin{eqnarray*}
\frac{\mu_{2,0}(C(0)\cap\mathbb{R}^2)\mu_{2,0}(C(0)\cap K_1\cap K_2)}{\mu_{2,0}(C(0)\cap K_1)\mu_{2,0}(C(0)\cap K_2)}\geq \frac{\mu_{2,0}(\mathbb{R}^2)\mu_{2,0}(M_{1}\cap M_{2,0})}{\mu_{2,0}(M_{1})\mu_{2,0}(M_{2,0})}. 
\end{eqnarray*}
This of course can be given by applying Theorem \ref{con}.

For $\beta>0$, Let $M^{\beta}_1$ be a width-decreasing set such that $Max(M^{\beta}_1)=2$ and such that it has the same axis of symmetry as the set $M_1$. This set is obtained by proceeding the symmetrisation procedure of the set $C(\beta)\cap K_1$ with respect to the measure $\mu_{2,0}$. Define the constant $c(\beta)>0$ such that :
\begin{eqnarray*}
\frac{\mu_{2,\beta}(C(\beta)\cap K_1)}{\mu_{2,\beta}(C(\beta)\cap \mathbb{R}^2)}=\frac{c(\beta)\mu_{2,\beta}(M^{\beta}_1)}{\mu_{2,\beta}(\mathbb{R}^2)},
\end{eqnarray*}
and use the symmetrisation procedure of Theorem \ref{con} to obtain a width-decreasing set $M_{2,\beta}$ such that $Max(M_{2,\beta})=2$, satisfying :
\begin{eqnarray*}
\frac{\mu_{2,\beta}(C(\beta)\cap K_2)}{\mu_{2,\beta}(C(\beta)\cap \mathbb{R}^2)}=\frac{\mu_{2,\beta}(M_{2,\beta})}{\mu_{2,\beta}(\mathbb{R}^2)},
\end{eqnarray*}

and

\begin{eqnarray*}
\frac{\mu_{2,\beta}(C(\beta)\cap K_1\cap K_2)}{\mu_{2,\beta}(C(\beta)\cap\mathbb{R}^2)}\geq \frac{c(\beta)\mu_{2,\beta}(M^{\beta}_1\cap M_{2,\beta})}{\mu_{2,\beta}(\mathbb{R}^2)}.
\end{eqnarray*}
By continuity, the set $M_{2,\beta}$ is constructed from $M_{2,0}$ by (slightly) moving it. 

Hence :
\begin{equation} \label{eqn:lthihi}
\frac{\mu_{2,\beta}(C(\beta)\cap\mathbb{R}^2)\mu_{2,\beta}(C(\beta)\cap K_1\cap K_2)}{\mu_{2,\beta}(C(\beta)\cap K_1)\mu_{2,\beta}(C(\beta)\cap K_2)}\geq \frac{c(\beta)\mu_{2,\beta}(\mathbb{R}^2)\mu_{2,\beta}(M^{\beta}_1\cap M_{2,\beta})}{c(\beta)\mu_{2,\beta}(M^{\beta}_1)\mu_{2,\beta}(M_{2,\beta})}.
\end{equation}

According to equation (\ref{eqn:lthihi}), there exists a $\beta_0$ such that (with respect to the cartesian coordinates associated to the measure $\mu_{2,\beta_0}$) the set $M^{\beta_0}_1$ is a \emph{nice convex set}.

Therefore, applying Theorem \ref{twodim} for this $\beta_0$, the proof of the Corollary \ref{corrr} follows.

%Let $\theta=0$ and for $i=1,2$, denote by $y_i$ the end point of the unit vector of the axis of the symmetric strips $X_i$ such that :
%\begin{eqnarray*}
%\frac{\mu_2(K_1\cap K_2)}{\mu_2(K_1)\mu_2(K_2)}\geq \frac{\mu_2(X_1\cap X_2)}{\mu_2(X_1)\mu_2(X_2)}.
%\end{eqnarray*}
%Without loss of generality we assume the $y_1\in [\pi/2,\pi]$. According to the proof of Theorem \ref{con}, for every cone $C(\theta)$ we can find symmetric strips $Y_1$ and $Y_2$ such that :
%\begin{eqnarray*}
%\frac{\mu_2(C(\theta)\cap K_1\cap K_2)\mu_2(C(\theta)}{\mu_2(K_1\cap C(\theta))\mu_2(K_2\cap C(\theta))}\geq \frac{\mu_2(\mathbb{R}^2)\mu_2(Y_1\cap Y_2)}{\mu_2(Y_1)\mu_2(Y_2)},
%\end{eqnarray*}
%where the end point of the unit vector of the axis of $Y_1$ coincides with $y_1$.

%Therefore, the proof of the Corollary follows by applying Lemma \ref{rot} followed by Theorem \ref{con}.
\begin{flushright}
$\Box$
\end{flushright}

\section{Proof of Theorem \ref{strong}}

We are ready to finalise the proof of theorem \ref{strong}. To remind the setting, $K_1, K_2$ are two symmetric convex bodies in $\mathbb{R}^n$ and $\mu_n=g(r)dx$ is a radial probability measure
.

For $i=1,2$ define the functions $f_i:u\in\mathbb{S}^{n-1}\to\mathbb{R}$ by
\begin{eqnarray*}
f_i(u)=\int_{0}^{x_i(u)}g(t)t^{n-1}dt
\end{eqnarray*}
where $x_i(u)$ is the length of the segment issuing from the origin in the direction $u$ where the other end touches the boundary of $K_i$.
Define the function $f_3:u\in\mathbb{S}^{n-1}\to\mathbb{R}$ by
\begin{eqnarray*}
f_3(u)=\int_{0}^{min\{x_1(u),x_2(u)\}}g(t)t^{n-1}dt.
\end{eqnarray*}
and finally, the constant function $f_4:u\in\mathbb{S}^{n-1}\to\mathbb{R}$ by
\begin{eqnarray*}
f_4(u)=\int_{0}^{+\infty}g(t)t^{n-1}dt.
\end{eqnarray*}
It is clear that for every $u\in\mathbb{S}^{n-1}$ we have 
\begin{eqnarray*}
f_4(u)\,f_3(u)\geq f_1(u)\,f_2(u).
\end{eqnarray*}

We proceed by contradiction: assume inequality (\ref{eqn:do}) holds for every pair of strips $S_1,S_2$ where the axis of at least one coincides with the $x$-axis and assume there exists $K_1, K_2$ such that
\begin{equation} \label{eqn:nab}
\mu_n(K_1\cap K_2)<\mu_n(K_1)\mu_n(K_2).
\end{equation}

According to equation (\ref{eqn:nab}), we can find a $C>0$ such that:
\begin{eqnarray*}
\frac{\displaystyle\int_{\mathbb{S}^{n-1}}f_1(u)d\sigma(u)}{\displaystyle\int_{\mathbb{S}^{n-1}}f_4(u)d\sigma(u)}> C> \frac{\displaystyle\int_{\mathbb{S}^{n-1}}f_3(u)d\sigma(u)}{\displaystyle\int_{\mathbb{S}^{n-1}}f_2(u)d\sigma(u)}.
\end{eqnarray*}
Define the two functions:
\begin{eqnarray*}
G_1(u)=f_1(u)-Cf_4(u),
\end{eqnarray*}
and
\begin{eqnarray*}
G_2(u)=Cf_2(u)-f_3(u).
\end{eqnarray*}
$G_1$ and $G_2$ are two continuous functions on $\mathbb{S}^{n-1}$ such that for $i=1,2$ we have:
\begin{eqnarray*}
\int_{\mathbb{S}^{n-1}}G_i(u)d\sigma(u)> 0.
\end{eqnarray*}
We now apply Lemma \ref{genlova} which provides us with a convex partition of $\mathbb{S}^{n-1}$ into geodesic segments and a family of $\sin^{n-2}$-affine probability measures on every geodesic segment of this partition for which we have:

\begin{equation} \label{eqn:tu}
\int_{\sigma}G_i(t)d\nu(t)> 0.
\end{equation}

Let $I$ be a geodesic segment of the partition given by Lemma \ref{genlova}. Take the two-dimensional cone $C(I)$ which contains the origin of $\mathbb{R}^n$ and the geodesic segment $I$. According to equation (\ref{eqn:tu}), for every geodesic segment $I$ of the partition we have:
\begin{equation} \label{eqn:neg}
\frac{\mu_2(C(I)\cap\mathbb{R}^2)\mu_2(C(I)\cap (K_1\cap K_2))}{\mu_2(C(I)\cap K_1)\mu_2(C(I)\cap K_2)}< 1.
\end{equation}

By applying Theorem \ref{con} for every $I$ geodesic segment of the partition given by Lemma \ref{genlova}, we find symmetric strips $S_1(I)$ and $S_2(I)$ such that:
\begin{eqnarray*}
\frac{\mu_2(C(I)\cap\mathbb{R}^2)\mu_2(C(I)\cap (K_1\cap K_2))}{\mu_2(C(I)\cap K_1)\mu_2(C(I)\cap K_2)}\geq \frac{\mu_2(\mathbb{R}^2)\mu_2(S_1(I)\cap S_2(I))}{\mu_2( S_1(I))\mu_2(S_2(I))}.
\end{eqnarray*}

% We are set to apply Theorem \ref{twodim} and its corollary \ref{con}. The symmetrisation which are used in this theorem and its corollary procures us with two symmetric strips $S_1$ and $S_2$ such that:
%\begin{eqnarray*}
%\frac{(\int_{I}f_4(t)d\nu(t))(\int_{I}f_3(t)d\nu(t))}{(\int_{I}f_1(t)d\nu(t))(\int_{I}f_2(t)d\nu(t))}\geq \frac{\mu_2(\mathbb{R}^2)\mu_2(S_1\cap S_2)}{\mu_2(S_1)\mu_2(S_2)}.
%\end{eqnarray*}
%If we can find one geodesic segment $I$ for which after the symmetrisation used in theorem \ref{twodim} and corollary \ref{con} we find two strips $S_1$ and $S_2$ with either $S_1$ or $S_2$ having $x$ as its axis, then we find a contradiction and the proof of this Theorem follows immediately. So we suppose that this is not the case and that for every $I$ geodesic of the partition, after the symmetrisation we find $S_1$ and $S_2$ none of them having $x$ as its axis.

If the axis of either $S_1(I)$ or $S_2(I)$ coincided with the $x$-axis, equation (\ref{eqn:neg}) and the assumption of the theorem would lead to a contradiction and the proof would be done. However, for the moment, we know nothing about the axis of the strips $S_i(I)$. To have a control on the axis of at least one strip, we need to work a bit harder.

%The issue at first glance is that nothing tells us that we may \emph{always} find geodesic segment $I$ in this partition such that the axis of either $S_1(I)$ or $S_2(I)$ coincides with the $x$-axis. Simple examples of this fact can be constructed- for example take for $K_1$ a cylinder in $\mathbb{R}^n$ and for $K_2$ a strip of widths equal to $\varepsilon$ and construct a simple partition of $\mathbb{S}^{n-1}$ by maximal spherical needles such that none of the strips in the cone over the spherical needles have an axis coinciding with the $x$-axis.

Note that by definition of $G_1$ and $G_2$ and due to the central symmetry of $K_1$ and $K_2$, on every half-sphere i.e. $\mathbb{S}^{n-1}_{+}$, we have:
\begin{equation} \label{eqn:have}
\int_{\mathbb{S}^{n-1}_{+}}G_i(u)d\sigma(u)=\frac{1}{2}\int_{\mathbb{S}^{n-1}}G_i(u)d\sigma(u).
\end{equation}
This means at first, we can freely cut the sphere in two in order to find the convex partition given by Lemma \ref{genlova}. In other words, the partition obtained by applying Lemma \ref{genlova} is not \emph{unique} and depends on the choice of the first cut of the sphere into two (we knew this already from the proof of Lemma \ref{genlova}). 

According to corollary \ref{corrr}, for every $\mathbb{S}^1\subset \mathbb{S}^n$, there exists a spherical needle (which will be denoted by $(J_\beta,\mu_{2,\beta})$ where $J_\beta\subset \mathbb{S}^1$ such that inequality (\ref{eqn:shahkar}) holds. We call such a needle a \emph{good needle}. Then we have:

\begin{lem} \label{close}
Let $n\geq 4$. For every $\varepsilon>0$, there exists a partition $\Pi$, a $J\in \Pi$, a \emph{good} needle $(J_\beta,\mu_{2,\beta})$ such that:
\begin{eqnarray*}
d(J,J_{\beta})\leq \varepsilon,
\end{eqnarray*}
where $d(J,J_{\beta})$ is the distance between the two needles defined in section $4$. The partitions are those partitions satisfying the assumption of Lemma \ref{genlova}.
\end{lem}
\emph{Proof of Lemma \ref{close}} :

%Choose a hemi-sphere $x^{\vee}_1$ in $\mathbb{S}^{n-1}$ and consider it as the first cut of the sphere $\mathbb{S}^{n-1}$. Let the center of this hemi-sphere be denoted by $x_1$. By ``a maximal spherical needle centered at $x_1$'', we mean a maximal geodesic containing $x_1$, such that the maximum point of the density function of its associated probability measure is at the point $x_1$. Let $I$ be a maximal spherical needle centered at $x_1$. Let $\beta>0$ and rotate the center $x_1$ by $\beta$ along the maximal geodesic segment $I$ (Choose a clock-wise rotation). Denote the new maximal spherical needle centered at $x_1(\beta)$ by $(I_{\beta},\nu_{\beta})$. Denote by $x^{\vee}_{1}(\beta)$ the hemi-sphere centered at $x_1(\beta)$.

We proceed by contradiction: there exists a $\varepsilon_0$ such that for every partition $\Pi$ obtained by applying Lemma \ref{genlova}, every needle $(I_1,\nu_1)$ in $\Pi$, every good needle $(I_2,\nu_2)$ we have 
\begin{equation} \label{eqn:teuf}
d((I_1,\nu_1),(I_2,\nu_2))\geq \varepsilon_0.
\end{equation}
Choose $(I_1,\nu_1)$ and $(I_2,\nu_2)$ satisfying (\ref{eqn:teuf}) and let $d((I_1,\nu_1),(I_2,\nu_2))=\varepsilon$. Furthermore, assume $\varepsilon$ is minimal when the pair $((I_1,\nu_1),(I_2,\nu_2))$ is chosen as above. A partition $\Pi$ exists (by assumption) such that the needle $(I_1,\nu_1)\in \Pi$. 

Let $\delta>0$ be such that $\delta<\varepsilon_0/2$. We choose $K_1$ (resp $K_2$) a $(1,\delta)$- constructing pancake for $(I_1,\nu_1)$ (resp $(I_2,\nu_2)$) such that $I_1\subset K_1\subset I_1+\delta$ and $I_2\subset K_2\subset I_2+\delta$. Therefore by definition, the Hausdorff distance between $K_1$ and $K_2$ is larger than $\varepsilon_0$. At the same time, by the minimality of $\varepsilon$, the intersection of $K_1$ and $K_2$ is non-empty. Let $w=I_1\cap I_2$, the point of intersection of $I_1$ and $I_2$.

By definition, the partition $\Pi$ is a (weak) limit of a finite partition $\Pi^N$, where $K_1\in \Pi^N$ is a constructing pancake. We now apply the algorithmic cutting procedure used in the proof of Lemma \ref{genlova} (by the aim of the Borsuk-Ulam Theorem) to cut the constructing pancakes of partition $\Pi^N$ in an appropriate manner to reach a \emph{new} partition different from $\Pi$ but still satisfying the properties of Lemma \ref{genlova}. Since $n\geq 4$, (we can) choose the cutting hemi-sphere to be a hemi-sphere with its boundary being a $(n-2)$-dimensional sphere, which is orthogonal to the geodesic segment $I_1$. We start to cut (solution of equation (\ref{eqn:borsuk}))
 with respect to this cutting hemi-sphere. After the first cut, we choose a new set $K^1_1$ which is the intersection of $K_1$ with the cutting hemi-sphere containing the point $w$. We carry on the cutting process each time with respect to a cutting-hemi-sphere, such that its boundary is orthogonal to $I_1$. At each stage we choose the hemi-sphere containing the point $w$. We denote the set obtained after cutting $j$ times by $K^j_1$. By the choice of the cutting hemi-spheres, the distance between the geodesic $I_2$ to the set $K^j_1$ is non-increasing as $j$ increases. Hence eventually, there will be a $K^j_1$ pancake such that its distance to a $J_\beta\subset I_2$ would be less than $\delta$. Moreover, according to the claim in the proof of Theorem \ref{con}, the needle $(J_\beta,\nu_{\beta})$, where $\nu_{\beta}$ is the probability measure obtained by restricting $\nu_2$ on $J_{\beta}$, is a \emph{good needle}. This is a contradiction with the definition of $\varepsilon_0$. The proof of Lemma \ref{close} is completed.

\begin{flushright}
$\Box$
\end{flushright}

Define the following function on the space of convexly-derived measures $\mathcal{MC}^{\leq 1}$, 
\begin{eqnarray*}
F:\mathcal{MC}^{\leq 1}\to\mathbb{R}_{+}
\end{eqnarray*}
by
\begin{eqnarray*}
F((I,\nu))=\frac{\mu_{2,\nu}(C(I)\cap \mathbb{R}^2)\mu_{2,\nu}(C(I)\cap K_1\cap K_2)}{\mu_{2,\nu}(C(I)\cap K_1)\mu_{2,\nu}(C(I)\cap K_2)},
\end{eqnarray*}
where $\mu_{2,\nu}$ is the measure written with respect to $\nu$ on the cone $C(I)$.

The function $F$ is a continuous function. For every $\varepsilon>0$, by compacity of the space of convexly derived measures $\mathcal{MC}^{\leq 1}$, there exists a $C(\varepsilon)>0$ such that for every spherical needle $(I_1,\nu_1)$ which is $\varepsilon$-close to a spherical needle $(I_0,\nu_0)$, and such that $F((I_0,\nu_0))\geq F((I_1,\nu_1))$ we have :
\begin{equation} \label{eqn:ineqq}
\frac{\mu_{2,\nu_1}(C(I_1)\cap\mathbb{R}^2)\mu_{2,\nu_1}(C(I_1)\cap K_1\cap K_2)}{\mu_{2,\nu_1}(C(I_1)\cap K_1)\mu_{2,\nu_1}(C(I_1)\cap K_2)}\geq \frac{C(\varepsilon)\mu_{2,\nu_0}(C(I_0)\cap\mathbb{R}^2)\mu_{2,\nu_0}(C(I_0)\cap K_1\cap K_2)}{\mu_{2,\nu_0}(C(I_0)\cap K_1)\mu_{2,\nu_0}(C(I_0)\cap K_2)}.
\end{equation}

For every $\varepsilon>0$, by applying Lemma \ref{close} combined with Corollary \ref{corrr}, we know that there exists a spherical needle $(I_1,\nu_1)$ which is $\varepsilon$-close to a spherical needle $(I_{\beta},\nu_{\beta})$ such that for the spherical needle $(I_{\beta},\nu_{\beta})$ we have :
\begin{eqnarray*}
F((I_{\beta},\nu_{\beta})\geq \frac{\mu_{2,\beta}(\mathbb{R}^2)\mu_{2,\beta}(S_1\cap S_2)}{\mu_{2,\beta}(S_1)\mu_{2,\beta}(S_2)},
\end{eqnarray*}
where $S_1$ and $S_2$ are symmetric strips in $\mathbb{R}^2$ such that at least one of their axis coincides with the $x(\beta)$-axis when writing down the measure $\mu_{2,\beta}$ in the Cartesian coordinates associated to this measure.

If we have $F((I_1,\nu_1))\geq F((I_{\beta},\nu_{\beta})$ then we will have :
\begin{eqnarray*}
F((I_1,\nu_1))&\geq& F((I_{\beta},\nu_{\beta})) \\
              &\geq& \frac{\mu_{2,\beta}(\mathbb{R}^2)\mu_{2,\beta}(S_1\cap S_2)}{\mu_{2,\beta}(S_1)\mu_{2,\beta}(S_2)}.
\end{eqnarray*}
And, by assumption on the symmetric strips (when at least one has an axis coinciding with the $x(\beta)$-axis) we shall have :
\begin{eqnarray*}
F((I_1,\nu_1))&\geq& \frac{\mu_{2,\beta}(\mathbb{R}^2)\mu_{2,\beta}(S_1\cap S_2)}{\mu_{2,\beta}(S_1)\mu_{2,\beta}(S_2)}\\
              &\geq& 1.
\end{eqnarray*}
Which is a contradiction as we have equation (\ref{eqn:neg}).

Therefore, we can assume that for every $\varepsilon>0$ and for every $(I_1,\nu_1)$ and $(I_{\beta},\nu_{\beta})$ as above, we have :
\begin{eqnarray*}
F((I_1,\nu_1))&\geq& C(\varepsilon)F((I_{\beta},\nu_{\beta}))\\
              &=& \frac{C(\varepsilon)\mu_{2,\beta}(C(I_{\beta})\cap\mathbb{R}^2)\mu_{2,\beta}(C(I_{\beta})\cap K_1\cap K_2)}{\mu_{2,\beta}(C(I_{\beta})\cap K_1)\mu_{2,\beta}(C(I_{\beta})\cap K_2)}\\
              &\geq& \frac{C(\varepsilon))\mu_{2,\beta}(\mathbb{R}^2)\mu_{2,\beta}(S_1\cap S_2)}{\mu_{2,\beta}(S_1)\mu_{2,\beta}(S_2)}\\
              &\geq& C(\varepsilon).
\end{eqnarray*}

Since $\lim_{\varepsilon\to 0}C(\varepsilon)= 1$, choosing $\varepsilon$ small enough and a spherical needle $J_{\beta}$ sufficiently close to $I_{\beta}$, the sets $C(J_{\beta})\cap K_i$ are sufficiently close to the sets $C(I_{\beta})\cap K_i$ . Therefore combining this with equation (\ref{eqn:neg}), we obtain:
\begin{eqnarray*}
1 &>&\frac{\mu_2(C(J_{\beta})\cap\mathbb{R}^2)\mu_2(C(J_{\beta})\cap (K_1\cap K_2))}{\mu_2(C(J_{\beta})\cap K_1)\mu_2(C(J_{\beta})\cap K_2)} \\
  &\geq& \frac{\mu_2(\mathbb{R}^2)\mu_2(S_1\cap S_2)}{\mu_2(S_1)\mu_2(S_2)},
\end{eqnarray*}
 
where the axis of either $S_1$ or $S_2$ coincides with the $x(\beta)$-axis.

By assumption of the theorem, for every two symmetric strips (such that the axis of either $S_1$ or $S_2$ coincides with the $x(\beta)$-axis), we have:
\begin{eqnarray*}
\frac{\mu_2(\mathbb{R}^2)\mu_2(S_1\cap S_2)}{\mu_2(S_1)\mu_2(S_2)}\geq 1.
\end{eqnarray*}
This provides us with a contradiction and the proof of Theorem \ref{strong} is completed.

\begin{flushright}
$\Box$
\end{flushright}

\section{A proof of the Gaussian Correlation Conjecture}
  
Recall the Gaussian Correlation Conjecture (or the Gaussian Correlation Problem as the conjecture was settled by T.Royen in \cite{roy}).

\begin{theo} \label{gcc}
Let $\gamma_n$ be the standard Gaussian measure defined by equation (\ref{eqn:gaus}). For every pair of symmetric convex bodies $K_1,K_2$, we have:
\begin{eqnarray*}
\gamma_n(K_1\cap K_2)\geq \gamma_n(K_1)\gamma_n(K_2).
\end{eqnarray*}
\end{theo}

\emph{Proof of Theorem \ref{gcc}}:

Apply theorem \ref{strong}. In order to prove Theorem \ref{gcc} for $n\geq 4$, it simply remains to demonstrate the following:
\begin{lem} \label{end}
For every two symmetric strips $S_1$ and $S_2$ in $\mathbb{R}^2$, where the axis of either $S_1$ or $S_2$ coincides with the $y$-axis, we have:
\begin{eqnarray*}
\mu_2(\mathbb{R}^2)\mu_2(S_1\cap S_2)\geq \mu_2(S_1)\mu_2(S_2).
\end{eqnarray*}
Where $\mu_2=C(n)\vert x\vert^{n-2}e^{-(x^2+y^2)/2}dx\,dy$ and where
\begin{eqnarray*}
C(n)=(\int_{-\pi/2}^{\pi/2}\cos(t)^{n-2}dt)^{-1}.
\end{eqnarray*}
\end{lem}

\emph{Proof of Lemma \ref{end}}:

At first, assume $S_1$ and $S_2$ to be orthogonal symmetric strips and their axes coincide with the axis of the Cartesian plane. In this case, $a,b\in\mathbb{R}_{+}$ exist, such that:
\begin{eqnarray*}
\mu_2(S_1)&=&C(n)\int_{-\infty}^{+\infty}\int_{-a}^{+a}\vert x\vert^{n-2}e^{-(x^2+y^2)/2}dx\,dy\\
          &=&C(n)\left(\int_{-\infty}^{+\infty}e^{-y^2/2}dy\right)\left(\int_{-a}^{+a}\vert x\vert^{n-2}e^{-x^2/2}dx\right).
\end{eqnarray*}
And
\begin{eqnarray*}
\mu_2(S_2)&=&C(n)\int_{-\infty}^{+\infty}\int_{-b}^{+b}\vert x\vert^{n-2}e^{-(x^2+y^2)/2}dy\,dx\\
          &=&C(n)\left(\int_{-b}^{+b}e^{-y^2/2}dy\right)\left(\int_{-\infty}^{+\infty}\vert x\vert^{n-2}e^{-x^2/2}dx\right).
\end{eqnarray*}          
Therefore:
\begin{eqnarray*}
\mu_2(S_1)\mu_2(S_2)&=&C(n)^2\left(\int_{-\infty}^{+\infty}e^{-y^2/2}dy\right)\left(\int_{-a}^{+a}\vert x\vert^{n-2}e^{-x^2/2}dx\right)\\
&&\left(\int_{-b}^{+b}e^{-y^2/2}dy\right)\left(\int_{-\infty}^{+\infty}\vert x\vert^{n-2}e^{-x^2/2}dx\right) \\
                     &=&C(n)^2\left(\int_{-\infty}^{+\infty}e^{-y^2/2}dy\right)\left(\int_{-\infty}^{+\infty}\vert x\vert^{n-2}e^{-x^2/2}dx\right)\\
                     &&\left(\int_{-a}^{+a}\vert x\vert^{n-2}e^{-x^2/2}dx\right)
                     \left(\int_{-b}^{+b}e^{-y^2/2}dy\right)\\
                     &=&\left(C(n)\int_{-\infty}^{+\infty}\int_{-\infty}^{+\infty}\vert x\vert^{n-2}e^{-(x^2+y^2)/2}dx\,dy\right)\\
                     &&\left(C(n)\int_{-b}^{+b}\int_{-a}^{+a}\vert x\vert^{n-2}e^{-(x^2+y^2)/2}dx\,dy\right)\\
                     &=&\mu_2(\mathbb{R}^2)\mu_2(S_1\cap S_2).
\end{eqnarray*}
Hence, in the very special case of orthogonal symmetric strips parallel to coordinate axes, we have equality in Lemma \ref{end}.

Now we switch to the general case: assume one of the symmetric strips, $S_2$, has $y$-as its axis and the other strip is arbitrary.

There exists $a,b,\rho$ such that
\begin{eqnarray*}
\frac{1}{2}\mu_2(S_2)=C(n)\int_{0}^{+\infty}\int_{-a}^{+a}\vert x\vert^{n-2}e^{-(x^2+y^2)/2}dx\,dy,
\end{eqnarray*}
and
\begin{eqnarray*}
\frac{1}{2}\mu_2(S_1)=C(n)\displaystyle\int_{0}^{+\infty}\displaystyle\int_{-b+\tan(\rho)x}^{+b+\tan(\rho)x}\vert x\vert^{n-2}e^{-(x^2+y^2)/2}dy\,dx.
\end{eqnarray*}
Define the function:
\begin{eqnarray*}
k(x)=\frac{\vert x\vert^{n-2}e^{-x^2/2}\displaystyle\int_{-b+\tan(\rho)x}^{+b+\tan(\rho)x}e^{-y^2/2}dy}{\vert x\vert^{n-2}e^{-x^2/2}\displaystyle\int_{-\infty}^{+\infty}e^{-y^2/2}dy}
\end{eqnarray*}
This function is clearly decreasing on $[0,+\infty)$, due to the fact that  $e^{-y^2}$ is a decreasing function.
Remember the following:
\begin{lem}[Gromov] \label{bg}
Let $f,g$ be two positive functions defined on $[0,\infty)$. If $f/g$ is non-increasing, then the function:
\begin{eqnarray*}
\frac{\displaystyle\int_{0}^{x}f(t)dt}{\displaystyle\int_{0}^{x}g(t)dt}
\end{eqnarray*}
is non-increasing on $[0,\infty)$.
\end{lem}
Therefore, we can apply Lemma \ref{bg} and conclude that the function:
\begin{eqnarray*}
\frac{\displaystyle\int_{0}^{X}\left(\vert x\vert^{n-2}e^{-x^2/2}\displaystyle\int_{-b+\tan(\rho)x}^{+b+\tan(\rho)x}e^{-y^2/2}dy\right)dx}{\displaystyle\int_{0}^{X}\left(\vert x\vert^{n-2}e^{-x^2/2}\displaystyle\int_{-\infty}^{+\infty}e^{-y^2/2}dy\right)dx}
\end{eqnarray*}
is a non-increasing function. The consequence of this fact is that:
\begin{eqnarray*}
\frac{\mu_2(S_1\cap S_2)}{\mu_2(S_1)}&=&
\frac{\displaystyle\int_{0}^{a}\left(\vert x\vert^{n-2}e^{-x^2/2}\displaystyle\int_{-b+\tan(\rho)x}^{+b+\tan(\rho)x}e^{-y^2/2}dy\right)dx}{\displaystyle\int_{0}^{+\infty}\left(\vert x\vert^{n-2}e^{-x^2/2}\displaystyle\int_{-b+\tan(\rho)x}^{+b+\tan(\rho)x}e^{-y^2/2}dy\right)dx}\\
                                    &\geq& \frac{\displaystyle\int_{0}^{a}\left(\vert x\vert^{n-2}e^{-x^2/2}\displaystyle\int_{-\infty}^{+\infty}e^{-y^2/2}dy\right)dx}{\displaystyle\int_{0}^{+\infty}\left(\vert x\vert^{n-2}e^{-x^2/2}\displaystyle\int_{-\infty}^{+\infty}e^{-y^2/2}dy\right)dx}\\
                                    &=&\frac{\mu_2(S_2)}{\mu_2(\mathbb{R}^2)}.
\end{eqnarray*}

This ends the proof of Lemma \ref{end}.
\begin{flushright}
$\Box$
\end{flushright}

The following Lemma is a well known and straightforward result:
\begin{lem} \label{dim}
If there exists an integer $N>2$ such that the Gaussian Correlation Conjecture holds true for every $n\geq N$, then the conjecture holds true in every dimension.
\end{lem}

Theorem \ref{strong} combined with Lemmas \ref{end} and \ref{dim}, complete the \emph{proof of the Gaussian Correlation Conjecture}. 
\begin{flushright}
$\Box$
\end{flushright}

\emph{Remarks}:

The fact that the axis of either $S_1$ or $S_2$ coincides with the $y$-axis is \emph{crucial} to us.

Indeed if we try to solve a case in which one of the symmetric strips has an $x$-axis and the other is simply an arbitrary strip, we encounter serious difficulties.  In that case, in order to apply Lemma \ref{bg} in the same way as we did previously, we would have to define a function:
\begin{eqnarray*}
k(y)=\frac{e^{-y^2/2}\displaystyle\int_{-b+\cot(\rho)y}^{+b+\cot(\rho)y}\vert x\vert^{n-2}e^{-x^2/2}dx}{e^{-y^2/2}\displaystyle\int_{-\infty}^{+\infty}\vert x\vert^{n-2}e^{-x^2/2}dx}.
\end{eqnarray*}
 However, the function $k(y)$ is not necessarily non-increasing. The behaviour of this function depends on the number $b$, and if $b<b_{max}$, where $b_{max}$ is the point at which the function $\vert x\vert^{n-2}.e^{-x^2/2}$ attains its maximum, then $k(y)$ is certainly not non-increasing.

\section{Epilogue}

\subsection{A Recipe for the Correlation Property of Radial Probability Measures} 
Theorem \ref{strong} goes beyond an alternative proof for the Gaussian Correlation Conjecture. Indeed this Theorem provides us with a recipe which enables one to check the correlation property of every radial probability measure. There are other radial probability measures for which their correlation property can be interesting. An example of such a measure is the Cauchy measure defined by 
\begin{eqnarray*}
\nu_n=C(1+\vert x\vert^2)^{-\frac{n+1}{2}},
\end{eqnarray*}
where $C$ is the normalisation constant. In \cite{memcor}, some partial results on the correlation property of $\nu_n$ was presented, however proving that $\nu_n$ satisfies the correlation property was shown to be a far reaching problem. 

Applying Theorem \ref{strong} for $\nu_n$, what remains before proving that $\nu_n$ satisfies the correlation property is a counter-part to the lemma \ref{end} with regards to $\nu_n$. 

The above remarks can lead us to ask about the ``if and only if'' case in Theorem \ref{strong}.  
\subsection{The ``If and Only If'' Case in Theorem \ref{strong}}

Theorem \ref{strong} would be even more interesting if one could replace the \emph{if} by \emph{if and only if}.

This would mean, if for a radial probability measure $\mu_n$, there are two strips $S_1,S_2$ with the axis of one coinciding with the $x$-axis, inequality (\ref{eqn:do}) is not verified, then $\mu_n$ \emph{does not} satisfy the correlation property.

%Suppose the Gaussian Correlation Conjecture is true for every two centrally-symmetric convex bodies. Let $K_1,K_2\subset \mathbb{R}^n$, two centrally-symmetric convex bodies are given. Suppose a (measurable) \emph{nice} partition of the sphere $\mathbb{S}^{n-1}$ by geodesic segments is also be given. Since the Gaussian Correlation Conjecture is assumed to be true, by repeating the arguments in the proof of Lemma \ref{fondaa} we can conclude that a segment $I$-element of the partition, a $\sin^{n-2}$-probability measure $\nu$ defined on $I$ exists such that
%\begin{eqnarray*}
%\frac{(\displaystyle\int f_4d\nu)(\displaystyle\int f_3 d\nu)}{(\displaystyle\int f_1d\nu)(\displaystyle\int f_2d\nu)}\geq 1,
%\end{eqnarray*} 
%where $f_i$, $i=1,2,3,4$ are defined as previously. Take the $2$-dimensional plane containing the origin and containing the segment $I$. Denote this plane by $P$. Then for $i=1,2$, $P_i=P\cap K_i$ are centrally symmetric convex sets. According to Theorem \ref{twodim}, two strips $S_1$ and $S_2$ in $\mathbb{R}^2$ exist such that  
%\begin{eqnarray*}
%\frac{\mu_2(\mathbb{R}^2)\mu_2(P_1\cap P_2)}{\mu_2(P_1)\mu_2(P_2)}\geq \frac{\mu_2(\mathbb{R}^2)\mu_2(S_1\cap S_2)}{\mu_2(S_1)\mu_2(S_2)}.
%\end{eqnarray*}
In order to give a proof for the \emph{``if and only if''} case, one has to somehow show that for every configuration of strips $S_1$ and $S_2$, there are always two $n$- dimensional symmetric convex sets $K_1$ and $K_2$, which (if we apply the localisation and then symmetrisation of Theorem \ref{twodim}) lead to the strips $S_1$ and $S_2$.

This seems like an interesting inverse type problem:
\begin{conj}
The \emph{if} in Theorem \ref{strong} can be replaced by an \emph{if and only if}.
\end{conj} 

The story does not end with the radial Correlation problems. The spherical localisation technique developed in Section $3$ (and mainly Lemma \ref{genlova}) can be applied to other problems:

\subsection{A Few Words on Mahler Conjecture}

This section is taken from \cite{memmahl}. The Mahler Conjecture suggests that for $K$ a centrally-symmetric convex body in $\mathbb{R}^n$, we have:
\begin{eqnarray*}
vol_n(K)\,vol_n(\mathring{K})&\geq& vol_n(I^n)\,vol_n(\mathring{I}^n) \\
                            &=& \frac{4^n}{\Gamma(n+1)},
\end{eqnarray*}
where $\mathring{K}$ is the polar of the symmetric convex body $K$. Let us see how one could attack this problem using the localisation techniques discussed in this paper. Apply the following steps:
\begin{itemize}
\item Consider the symmetric convex sets $K$ such that for \emph{almost} all the sections $P\cap K$, where $P$ is a $2$-dimensional plane containing the origin we have:
\begin{eqnarray*}
 \mathring{(P\cap K)}=P\cap\mathring{K}.
\end{eqnarray*}

\item Study the two-dimensional variational problem, which is to minimise:
\begin{eqnarray*}
\mu_2(S)\,\mu_2(\mathring{S}),
\end{eqnarray*}
where $\mu_2=C(n)\vert x\vert^{n-2}dx\,dy$ in $\mathbb{R}^2$ over all two-dimensional symmetric convex sets $S$, where $C(n)$ is the appropriate normalisation constant. In order to answer this two-dimensional variational problem, one could consult \cite{fradmey} (pointed out by M.Fradelizi).
\item Denote the minimum obtained in the previous step by $g(n)$. Applying Lemma \ref{genlova} we directly obtain that for every symmetric convex set $K$ in $\mathbb{R}^n$, we have:
\begin{eqnarray*}
vol_n(K)\,vol_n(\mathring{K})&\geq& m(n) \\
                             &=&vol_{n-1}(\mathbb{S}^{n-1})^2\,g(n),
\end{eqnarray*}
where $vol_{n-1}(\mathbb{S}^{n-1})$ is the $(n-1)$-dimensional volume of the canonical sphere $\mathbb{S}^{n-1}$.
\item $m(n)$ is a lower bound for the Mahler volume of $K$. 

Compare $m(n)$ with $\frac{4^n}{\Gamma(n+1)}$ to see how far we are from the lower bound suggested by the Mahler Conjecture.
\end{itemize}

\bibliographystyle{plain}
\bibliography{gcor}

\end{document}